\documentclass[11pt]{article}
\usepackage{amsfonts}
\usepackage{amssymb}
\begin{document}
\newtheorem{lemma}{Lemma}[section]
\newtheorem{theorem}{Theorem}[section]
\newtheorem{proposition}{Proposition}[section]
\newtheorem{corollary}{Corollary}[section]
\newtheorem{example}{Example}[section]
\newtheorem{examples}{Examples}[section]
\def\dim{\mathop{\rm dim}\nolimits}
\def\deg{\mathop{\rm deg}\nolimits}
\def\codim{\mathop{\rm codim}\nolimits}
\def\char{\mathop{\rm char}\nolimits}
\def\rank{\mathop{\rm rank}\nolimits}
\def\rn{\mathop{\rm rk}\nolimits}
\def\id{\mathop{\rm id}\nolimits}
\def\rnk{\mathop{\overline{\rm rk}}\nolimits}
\def\dg{\mathop{\rm dg}\nolimits}
\def\ind{\mathop{\rm ind}\nolimits}
\def\Hom{\mathop{\rm Hom}\nolimits}
\def\Hilb{\mathop{\rm Hilb}\nolimits}
\def\Im{\mathop{\rm Im}\nolimits}
\def\GL{\mathop{\rm GL}\nolimits}
\def\cA{\mathop{\cal A{}}\nolimits}
\def\cB{\mathop{\cal B{}}\nolimits}
\def\cD{\mathop{\cal D{}}\nolimits}
\def\cV{\mathop{\cal V{}}\nolimits}
\def\cQ{\mathop{\cal Q{}}\nolimits}
\def\cP{\mathop{\cal P{}}\nolimits}
\def\cU{\mathop{\cal U{}}\nolimits}
\def\cO{\mathop{\cal O{}}\nolimits}
\def\cG{\mathop{\cal G{}}\nolimits}
\def\cF{\mathop{\cal F{}}\nolimits}
\def\cK{\mathop{\cal K{}}\nolimits}
\def\cT{\mathop{\cal T{}}\nolimits}
\def\cZ{\mathop{\cal Z{}}\nolimits}
\def\cN{\mathop{\cal N{}}\nolimits}
\def\cM{\mathop{\cal M{}}\nolimits}
\def\cH{\mathop{\cal H{}}\nolimits}
\def\cX{\mathop{\cal X{}}\nolimits}
\def\cY{\mathop{\cal Y{}}\nolimits}
\def\cW{\mathop{\cal W{}}\nolimits}
\def\cI{\mathop{\cal I{}}\nolimits}
\def\cR{\mathop{\cal R{}}\nolimits}
\def\cC{\mathop{\cal C{}}\nolimits}
\def\cE{\mathop{\cal E{}}\nolimits}
\def\cS{\mathop{\cal S{}}\nolimits}
\def\cSN{\mathop{\cal SN{}}\nolimits}
\def\cSE{\mathop{\cal SE{}}\nolimits}
\def\cSD{\mathop{\cal SD{}}\nolimits}
\def\cSC{\mathop{\cal SC{}}\nolimits}
\def\cSU{\mathop{\cal SU{}}\nolimits}
\def\pf{\mbox{\bf Proof \ }}
\def\pfA{\mbox{\bf Proof A \  }}
\def\pfB{\mbox{\bf Proof B \ }}
\newcommand{\N}{\mbox{$\mathbb N$}}
\newcommand{\C}{\mbox{$\mathbb C$}}
\newcommand{\R}{\mbox{$\mathbb P$}}
\newcommand{\Z}{\mbox{$\mathbb Z$}}
\newcommand{\Q}{\mbox{$\mathbb Q$}}
\newcommand{\n}{\mbox{${\displaystyle \frac{n}{2} }$}}
\newcommand{\m}{\mbox{${\displaystyle \Big[ \frac{n}{2} \Big] }$}}
\newcommand{\mm}{\mbox{$ \Big[ \frac{n}{2} \Big] $}}
\hyphenation{as-so-cia-ted con-si-de-ring sti-mu-la-ting}
\title{On the maximum nilpotent orbit which intersects the centralizer of a matrix}
\author{{Roberta Basili}
\\{\small{\em Liceo Annesso al Convitto Nazionale "Principe di Napoli", Assisi, Italy}}}
\date{}
\maketitle \noindent
\begin{abstract} We introduce a method to determine the maximum nilpotent orbit
which intersects a variety of
nilpotent matrices described by a strictly upper triangular matrix over a polynomial ring. We show that the result only depends on the ranks of its submatrices and we introduce conditions on a subvariety so that it intersects the same orbit. Then we describe a maximal nilpotent subalgebra of
the centralizer of any nilpotent matrix; the previous method
allows us to show that the maximum nilpotent
orbit which intersects that centralizer only depends on which entries
are identically $0$ in that subalgebra. The aim of the paper is to prove a simple algorithm for the determination of the maximum nilpotent orbit
which intersects that centralizer, which was conjectured by Polona Oblak.
\vspace{1mm}
\vspace{1mm} \newline {\em Keywords:} commuting
matrices, nilpotent orbits, maximal orbit.
\end{abstract}
\section{Introduction}\label{introduction} Let $
M(n,K) $ be  the affine space of all the $n\times n$ matrices over
an infinite field $K$ and let $ \mbox{\rm GL}\, (n,K) $ be the
open subset of $M(n,K)$ of all the $n\times n$ nonsingular
matrices. We will denote by $ N (n,K)$  the subvariety of $M(n,K)$
of all the $n\times n$  nilpotent matrices and by $\cN $ the
affine subspace of $M(n,K)$ of all the strictly upper triangular
matrices. To any element of $N(n,K)$ it corresponds a partition of
$n$; we fix a matrix $J \in \cN $ with Jordan canonical form and
denote by $\mu _1\geq \mu _2\geq \cdots \geq \mu _t $ the orders
of the Jordan blocks of $J$, hence $ B=(\mu _1,\ldots ,\mu _t)$ is
the partition of $n$ associated to $J$ (which can be identified
with the orbit of $J$ under the action of GL$(n,K)$).\newline If $
J' $ is another element of $N(n,K)$, we denote by $\mu '_1\geq
\cdots \geq \mu '_{t'}$ the orders of its Jordan block and we set
$ B'=(\mu '_1,\ldots ,\mu '_{t'})$. Then $B=B'$ if and only if
$\rank J^m=\rank (J')^m$ for all $m\in \N $. It is said that $
B<B'$ if $\, \rank J^m\leq \rank (J')^m\, $ for all $ m\in \N$ and
there exists $ m\in \N$ such that $\, \rank J^m\, <\, \rank
(J')^m$. The following claim is due to Hesselink (\cite{He},
1976): $B<B'\,  $ iff $  B$, as an orbit, is contained in the
closure of the orbit $ B'$. \newline For $i\in \N $ we set $\mu
_i=0$ if $i>t$ and $\mu'_i=0$ if $i>t'$; then
$$ B\leq B' \qquad \Longleftrightarrow \qquad {\displaystyle \sum_{i=1}^l\mu _i \leq
\sum_{i=1}^l\mu '_i}\quad \mbox{\rm for all } \ l\in \N
$$ where equality holds in the first relation iff it holds in the
second one for all $l\in \N $. \begin{examples}
 {\rm $ (6,4,3)<(6,5,2)<(6,6,1),$ \hspace{1mm} $(5,3,2,1)<(6,3,1,1)<(6,4,1)$;\hspace{1mm} $\,(6,5,2)$ and $(7,3,3)$
 cannot be compared, the same for $(6,5,4,3)$
and $(6,6,2,2,2)$.}\end{examples} We will denote by $ {\cal
C }_B$ the centralizer of $J$ and by $ \cN _B$ the algebraic
subvariety of $\cC _B$ of all the nilpotent matrices. We recall
the following result, whose proof is a consequence of Wedderburn's
Theorems.
\begin{lemma}\label{algebra} If $\cU $ is a finite dimensional
algebra over an infinite field then the subvariety
of all the nilpotent elements of $\cU $ is irreducible.\end{lemma} By Lemma \ref{algebra}
(see also for example Lemma 2.3 of \cite{Bas}) $\cN _B$ is
irreducible; for
$m\in \N $ the subvariety of $\cN _B$ of all $X$ such that $\rank
X^m$ is the maximum possible is open and, by the irreducibility of $\cN _B$, the intersection of these open subsets obtained for $m\in
\N $ is not empty.
 Hence there is a maximum partition for the elements of
$\cN _B$ and the subset of the elements which have this partition
is open (dense) in $\cN _B$.
 Then we can define a map $Q$ in the set of
the orbits of $n\times n$ nilpotent matrices (or partitions of
$n$) which associates to any orbit $B$ the maximum nilpotent orbit
which intersects $\cN _B$.\newline  For $s\in \N -\{ 0\} $ let $q$
and $r$ be the quotient and the remainder of the division of $n$
by $s$; then if $J$ is the $n\times n$ Jordan block we have that
$J^s$ has $r$ Jordan blocks of order $q+1$ and $s-r$ Jordan blocks
of order $q$. Hence the partition $B$ is almost rectangular (that
is $\mu _1-\mu _t\leq 1$) iff $J$ is conjugated to a power of the
$n\times n$ Jordan block. This implies that if $B$ is almost
rectangular then $Q(B)=(n)$. The converse of this claim is also
true; it is a consequence of the next Proposition which we are
going to explain.
There exist $p\in \N $ and almost rectangular
partitions $B_1,\ldots ,B_p$ of numbers less than or equal to $n$
such that $B=(B_1,\ldots ,B_p)$; we denote by $r_B$ the minimum of
all $p$ with this property (sometimes it can be obtained with
different choices of $B_1,\ldots ,B_{r_B}$). \begin{examples}
{\rm If $B=(5,4,3,1,1)$ we have $r_B=3$, which is
obtained with $B_1=(5)$, $B_2=(4,3)$, $B_3=(1,1)$ or with
$B_1=(5,4)$, $B_2=(3)$, $B_3=(1,1)$. If $B=(9,7,5,1)$ we have
$r_B=4$.} \end{examples} Let $s_B$ be the maximum value of $l$ for which there
exists
 a subset $\{ i_1,\ldots , i_l\} $ of $ \{ 1,\ldots ,t\} $ such that
 $i_1< \ldots <i_l$ and $\mu _{i_1} -\mu _{i_l} \leq 1$. The
 following Propositions were proved in {\cite{Bas} (Propositions 2.4
 and 3.5).
\begin{proposition} \label{R1} There exists a not empty open subset of $\cN
_B$ such that if $A$ belongs to it then $\mbox{\rm rank}\ A=n-r_B$
(that is $A$ has $r_B$ Jordan blocks).
\end{proposition}
\begin{proposition}\label{R2} We have that $\, \mbox{\rm rank }\;  (A^{s_B})^m\leq \mbox{\rm
rank } \; J^m\, $ for all $A\in \cN _B$ and $m\in \N
$.\end{proposition} Let $B=(B_1,\ldots ,B_{r_B})$
where $B_i$ is almost rectangular for $i=1,\ldots , r_B$; let
$n_i$ be the sum of the numbers of $B_i$ and let $\widetilde B
=(n_1,\ldots , n_{r_B})$. The following result is a consequence of
Proposition \ref{R2} (see Theorem 1.11 of  \cite{BasI}).
\begin{proposition}\label{RI} If $\, s_B=|B_i| $ for $ \
i=1,\ldots ,r_B\ $ then  $\ Q(B)=\widetilde B$.\end{proposition}
\begin{example} {\rm If $B=(5,4,4,2,2,1)$ we have
$\widetilde B=(13,5)$ and $Q(B)=\widetilde B$.}\end{example}
\begin{corollary}\label{corR} $\ Q(B)=B\ $ iff $\, r_B=t\, $, that is
 $ \,  \mu _{i}-\mu _{i+1}>1\, $  for
$ \, i=1,\ldots ,t-1\, $. \end{corollary}
 Let $\{ q_1,q_2,\ldots ,q_u\} $ be the
ordered subset of $\{ 0,\ldots ,t\} \ $  such that $\ q_u=t\ $ and
$$\mu _1=\mu _{q_1}\neq \mu _{q_{1}+1}=\mu _{q_{2}}\neq \mu _{q_2+1}=\cdots \neq
\mu _{q_{u-1}+1}= \mu _{q_u}\ $$  (for example if $B=(6,6,6, 6,
5,2, 2, 1)$ we have that   $q_1=4$, $q_2=5$, $q_3=7$, $q_4=8$). If
we set $q_0=0$ then $J$ has  $ q_i-q_{i-1}$ Jordan blocks of order
$ \mu _{q_i}$ for $i=1,\ldots ,u$.
 We
will write the partition $(\mu _1,\ldots ,\mu _t)$ also as $$ (\mu
_{q_1}^{q_1},\mu _{q_2}^{q_2-q_1}\ldots ,\mu _{q_{u}
}^{q_{u}-q_{u-1}}) \ .$$ We consider the subset of  $ \{ 1,\ldots
,u\} \times \{ 0,1\} $ of all $(i, \epsilon )$ such that $\mu
_{q_i}-\mu _{q_{i+\epsilon }}\leq 1$ (that is such that $\epsilon
=1$ iff $i<u$ and $\mu _{q_{i+1}}=\mu _{q_i}+1$); then we consider
the map from this subset to $\N $ defined by
$$(i, \epsilon ) \quad \longmapsto \quad 2q_{i-1}+\mu
_{q_i}(q_i-q_{i-1})+\epsilon \, \mu _{q_{i+1}}(q_{i+1}-q_{i})\ .$$
We denote by $\omega _1$ the maximum of the image of this map.
Polona Oblak in \cite{Obl} (2007) proved the following result.
\begin{theorem} \label{O} The maximum index of nilpotency for an element of $\cN
_B$, that is the first number of the partition $Q(B)$, is $\
\omega _1$.
\end{theorem} We will denote by $(\tilde i, \tilde {\epsilon })$ any preimage of $\omega _1$ with
respect to the previous map. \begin{example}
{\rm If $B=(5^2,4,3^4,2,1)$ we can only have $\tilde i=2$, $\tilde
{\epsilon }=1$ and the first number of $Q(B)$ is $\omega_1=2\times
2 + 4+3\times 4=20$.}\end{example} The canonical basis
$\Delta _B $ of $K^n$ will be written in the following way:
$$\Delta _B=\{ \ v_{\mu _{q_i},j}^{\mu _{q_i}}, v_{\mu _{q_i},j}^{\mu _{q_i}-1},\ldots
, v_{\mu _{q_i},j}^{1},\ i=1,\ldots ,u,\ j=q_i-q_{i-1},\ldots ,1
\} \ .$$ For example if $B=(5,3,3,2,1)=(5^1,3^2,2^1,1^1)$ we write
$$\Delta _B=\{ \ {v_{5,1}^5,v_{5,1}^4,\ldots
,v_{5,1}^1},{v_{3,2}^3,v_{3,2}^2, v_{3,2}^1},
{v_{3,1}^3, v_{3,1}^2, v_{3,1}^1},{
v_{2,1}^2, v_{2,1}^1},{ v_{1,1}^1}\ \} \
.\vspace{2mm}$$ Let $\Delta _B^{\circ }$ be the union of the
following subsets:
$$\Delta _B^{\circ ,1}=\{ v_{\mu _{q_i},j}^1\ |\ j=q_i-q_{i-1},\ldots ,1\ , \ i=1,\ldots
,\tilde i+\tilde{\epsilon}\} \ , $$
$$\Delta _B^{\circ ,2}=\{
v_{\mu _{q_i},j}^l \ | \ j=q_i-q_{i-1}, \ldots ,1 \ ,\ l=1,\ldots
,\mu _{q_i} \ ,\ i=\tilde i,\tilde i+\tilde {\epsilon }\} \ ,$$
$$\Delta _B^{\circ ,3}=\{ v_{\mu
_{q_i},j}^{\mu _{q_i}}\ |\ j=q_i-q_{i-1},\ldots ,1\ , \ i=1,\ldots
,\tilde i+\tilde{\epsilon} \} \ . $$ Then $\big| \Delta _B^{\circ
}\big| =\omega _1$.  Let $\widehat B$ be the partition obtained
from $B$ by cancelling the powers $\mu _{q_{i}}^{q_{i}-q_{i-1}}$
for $i=\tilde i,\tilde i+\tilde{\epsilon}$ and decreasing by $2$
the numbers $\ \mu _{q_i}$ for $i=1,\ldots ,\tilde i-1$, that is:
$$\widehat B= \big( (\mu _{q_1}-2)^{q_1}, \ldots ,(\mu
_{q_{\tilde i-1}}-2)^{q_{\tilde i-1}-q_{\tilde i-2}}, \mu
_{q_{\tilde i+\tilde{\epsilon}+1}}^{q_{\tilde
i+\tilde{\epsilon}+1}-q_{\tilde i+\tilde{\epsilon}}}, \ldots , \mu
_{q_u}^{q_u-q_{u-1}}\big)  \ .$$  Let $Q(B)=(\omega_1,\ldots
,\omega_z)$; let $Q(\widehat {B})=(\widehat {\omega } _1, \ldots
,\widehat{\omega }_{\hat z})$. One of the main aims of this paper is to
prove the
following Theorem.
\begin{theorem}\label{eend} The maximum partition $(\omega_1,\omega_2,\ldots ,\omega_z)$ which is associated to
 elements of $ {\cN }_B$ is  $(\omega_1, \widehat{\omega } _1, \widehat{\omega}_2,\ldots , \widehat
{\omega}_{\hat z})$, that is $Q(B)=(\omega_1, Q(\widehat{B}))$.
\end{theorem}
Theorem \ref{eend} leads to an algorithm for the determination of
the maximum partition which is associated to elements of $\cN _B$
for any partition $B$.  \begin{example} {\rm If $B=(15,13,5,4,3^2,2,1)$ we have $\tilde i=q_{\tilde i}=\mu _{q_{\tilde i}}=4$, $\omega_1=4+3\cdot 2+2\cdot 3=16$ and
$\widehat B=(15-2,13-2,5-2,2,1)=(13,11,3,2,1)$. Then $\widehat {\widehat B}=(11,3,2,1)$ and $\widehat {\widehat {\widehat B}}=(3,2,1)$. Since the maximum partition of the
elements of $\cN _{\widehat {\widehat {\widehat B}}}$ is $(5,1)$, we get that the maximum
partition of the elements of $\cN _B$ is $(16,13,11,5,1)$.}\end{example}
The conjecture expressed by Theorem \ref{eend} was communicated by its author
Polona Oblak during the meeting "Fifth
Linear Algebra Workshop" which was held in Kranjska Gora (May 27 - June 5, 2008).
After that, A. Iarrobino and L. Khatami wrote a paper on the
inequality $(\omega_1,\omega_2,\ldots ,\omega_z)\geq (\omega_1,
\widehat{\omega } _1, \widehat{\omega}_2,\ldots , \widehat
{\omega}_{\hat z})$ (see \cite{KhI}). L. Khatami in \cite{Kha} proved the uniqueness of the
result of the algorithms which follow from Theorem \ref{eend}
(since there can be different choices of $(\tilde i ,\tilde
{\epsilon })$), in \cite{Khat}
gave a formula for $\omega_z$.\newline In Section \ref{orbit}
we introduce a method for the determination of the
maximum nilpotent orbit which intersects a subvariety of  $\cN $ and we show the following theorem.
\begin{theorem}\label{equalorbit} If $\Upsilon $ and $\Upsilon '$ are two rational functions from the same affine variety
$\cA $ to $\cN $ such that any submatrix of $\Upsilon $ has the same rank as the corresponding submatrix of $\Upsilon '$ then the maximum nilpotent orbit which intersects $\Upsilon (\cA )$ is the same as the maximum nilpotent orbit which intersects $\Upsilon '(\cA )$.\end{theorem}
In Section \ref{centralizer} we introduce the variety $\cSN _B$, a maximal nilpotent subalgebra of $\cN _B$. We consider the minimal subspace $\cU _B$ of
$M(n,K)$ containing $\cC _B$ which is defined
by the condition that some coordinates are $0$ and we denote by $\cE _B$ the subvariety of all the nilpotent matrices of $\cU _B$. Besides describing $\cSN _B$ and the corresponding maximal nilpotent subalgebra $\cSE _B$ of $\cU _B$, by the results of Subsection \ref{P} (which focuses on a property which implies the hypothesis of Theorem \ref{equalorbit}) we prove the following theorem. \begin{theorem}\label{indToep} The maximum orbit intersecting
$\cSE _B$ ($\cE _B$) is $\, Q(B)$.
\end{theorem}
In Section \ref{conjecture} we give a proof of Theorem
\ref{eend} in which we replace $\cSN _B$ with $\cSE _B$, according to Theorem \ref{indToep}; nevertheless that proof could easily be rewritten with the same arguments and considering only the elements of $\cSN _B$. In Section \ref{fifth} we give
another proof of Theorem 3.2 of \cite{Kos},
which states the idempotency of $Q$,
besides recalling other results on this subject.
\section{On the maximum nilpotent orbit intersecting a subvariety of $\cN $}\label{orbit}
\subsection{The rational function $F(U)$}\label{rational}
We denote by $\Xi =(\xi _{i,j})$, $\, i,j=1,\ldots ,n$ the matrix of the
coordinates of the affine space $M(n,K)$. For any map $\Phi $ from $\{ 1,\ldots
,n\} $ to $\{ 2,\ldots ,n+1\} $ such that $\Phi (i)>i$ for $i=1,\ldots ,n$ let $\cN _{\Phi }$ be the quasi projective subvariety of $M(n,K)$ defined by the following conditions on the entries of $\Xi $:
$$ \begin{array}{l}\xi _{i,l}=0 \qquad \mbox{\rm if }\  \ l<\Phi (i)\ \qquad \mbox{\rm and } \qquad
  \xi _{i,\Phi (i)}\neq 0 \qquad  \mbox{\rm if } \ \ \Phi (i)\neq n+1
\ .\end{array}  $$
Similarly, for any map $\Psi $ from $\{ 1,\ldots
,n\} $ to $\{ 0,\ldots ,n-1\} $ such that $\Psi (i)<i$ for $i=1,\ldots ,n$ we can define $\cN ^{\Psi }$ as the quasi projective subvariety of $M(n,K)$ defined by the following conditions on the entries of $\Xi $:
$$ \begin{array}{l}\xi _{l,i}=0 \qquad \mbox{\rm if }\  \  l>\Psi (i)\ \qquad \mbox{\rm and  } \qquad
 \xi _{\Psi (i),i}\neq 0 \qquad  \mbox{\rm if } \ \ \Psi (i)\neq 0
\ .\end{array}  $$
Anyway, we will consider only the variety $\cN _{\Phi } $ (results similar to those in this section concerning $\cN _{\Phi } $ could be obtained for $\cN ^{\Psi } $).
The variety $\cN _{\Phi } $ is contained in $ \cN $ (the subvariety of all the $n\times n$
strictly upper triangular matrices); let $N_{\Phi }$ be the closure of
$\cN _{\Phi }$ in $\cN $
and let $\Xi ^{\Phi }$ be the matrix of the coordinates of $N_{\Phi }$, that is $\Xi ^{\Phi }=(\xi _{i,j})$, $i,j=1,\ldots ,n$, $\xi _{i,j}=0$ iff $j<\Phi (i)$.  Let $U=(u_{i,j})$,
$i,j=1,\ldots ,p$, $\ p>1$ be a submatrix of $\Xi ^{\Phi } $
with the following property:\begin{itemize} \item[$\star $)] $u_{i,k}=0$ if $k\in \{ 1,\ldots
,i-1\} $, if $p>2$ then $u_{i,i}\neq 0$ for $i=2,\ldots ,p-1$.
\end{itemize} If $h\in \{ 1,\ldots ,p-1\} $ and $k\in \{ h+1,
\ldots ,p\} $ let $U_{(h)}^{(k)}$ be the submatrix of $U$ obtained by
choosing the columns and the rows of $U$ of indices $h,\ldots ,k$.
The matrix $U_ {(h)}^{(k)}$ for $h=1,\ldots ,p-1$ and
$k=h+1,\ldots ,p$ has still the property $\star $). For
$h=1,\ldots ,p-1$ we set $U_{(h)}=U_{(h)}^{(p)}$; for $k=2,\ldots ,p$ we set
$U^{(k)}=U_{(1)}^{(k)}$. Obviously we have $U_{(1)}^{(p)}=U$.\newline For $h=1,\ldots
, p-1$ we define in the following way a rational function $F(U_{(h)})$
over $N_{\Phi } $ in the entries of $U_{(h)}$: we set $F(U_{(p-1)})=u_{p-1,p}$
and for $h=p-2,\ldots , 1$ we set
$$F(U_{(h)})=u_{h,p}-\sum_{k=h+1}^{p-1}u_{h,k}(u_{k,k})^{-1}F(U_{(k)})\
.$$ \begin{example}\label{ex1}{\rm \mbox{\rm If}\  \  $$U=\left(
\begin{array}{ccccc} u_{1,1} & u_{1,2} & u_{1,3} & u_{1,4} &
u_{1,5} \\ & u_{2,2} & u_{2,3} & u_{2,4} & u_{2,5} \\ & & u_{3,3}
& u_{3,4} & u_{3,5} \\ & & & u_{4,4} & u_{4,5} \\ & & & &
u_{5,5}\end{array}\right)\ \ \mbox{\rm where } u_{2,2},u_{3,3},
u_{4,4}\neq 0 \  \mbox{\rm then:}$$ $$F(U_{(4)})=u_{4,5}\ ,\quad
F(U_{(3)})=u_{3,5}-u_{3,4}(u_{4,4})^{-1}F(U_{(4)})=u_{3,5}-u_{3,4}(u_{4,4})^{-1}u_{4,5}\
,$$
$$F(U_{(2)})=u_{2,5}-u_{2,4}(u_{4,4})^{-1}F(U_{(4)})-u_{2,3}(u_{3,3})^{-1}F(U_{(3)})=$$
$$=u_{2,5}-u_{2,4}(u_{4,4})^{-1}u_{4,5}-u_{2,3}(u_{3,3})^{-1}(u_{3,5}-u_{3,4}(u_{4,4})^{-1}u_{4,5})\, ,$$
$$F(U)=u_{1,5}-u_{1,4}(u_{4,4})^{-1}F(U_{(4)})-u_{1,3}(u_{3,3})^{-1}F(U_{(3)})-u_{1,2}(u_{2,2})^{-1}F(U_{(2)})=$$ $$=u_{1,5}-
u_{1,4}(u_{4,4})^{-1}u_{4,5}-u_{1,3}(u_{3,3})^{-1}(u_{3,5}-u_{3,4}(u_{4,4})^{-1}u_{4,5})+ $$ $$-u_{1,2}(u_{2,2})^{-1}[
u_{2,5}-u_{2,4}(u_{4,4})^{-1}u_{4,5}-u_{2,3}(u_{3,3})^{-1}(u_{3,5}-u_{3,4}(u_{4,4})^{-1}u_{4,5})]\ .$$}\end{example}
The rational function $F(U)$ is a sum of fractional monomials, each one of degree $1$ with respect to the entries of the last column (of the first row) and of degree $0$ with respect to the entries of the other columns (rows). \newline
For $r\in \{ 1,\ldots ,p-1\} $ and $s\in \{ r+1,\ldots ,p\} $ we denote by $P_h(r,s)$ the coefficient of $u_{r,s}$ in $F(U_{(h)})$, for $h=1,\ldots ,p-1$. The following claim is obvious by
the definition of $F(U_{(h)})$. \begin{lemma}\label{l11}  $P_h(1,1)=P_h(p,p)=0$ for $h=1,\ldots ,p-1$; if $\; h\in \{ 2,\ldots ,p-1\} $ and $r\in \{ 1,\ldots ,h-1\} $, $s\in \{ r+1,\ldots ,p\} $ then $P_h(r,s)=0$.\end{lemma}
By the definition of $F(U_{(h)})$ we also get the following result.
\begin{lemma}\label{coeff} If $r\in \{ 1,\ldots ,p-1\} $ and $s\in \{ r+1,\ldots ,p\} $ then:
$$P_1(r,s)=\left\{ \begin{array}{cc} 1 & \mbox{\rm if } \ (r,s)=(1,p) \\ \\ -(u_{r,r})^{-1} F(U^{(r)}) & \mbox{\rm if $r\neq 1$ and $s=p$} \\ \\ -(u_{s,s})^{-1}F(U_{(s)}) & \mbox{\rm if $r=1$ and $s\neq p$} \\ \\ (u_{r,r})^{-1}(u_{s,s})^{-1}F(U_{(s)})F(U^{(r)}) & \mbox{\rm otherwise } \ . \end{array}\right. $$\end{lemma}
\pf By Lemma \ref{l11} the claim is true if $r=1$, hence we can prove it by induction on $r$. Let $r\neq 1$; we have that \begin{equation}\label{formaF} F(U)=u_{1,p}-\sum _{k=2}^{p-1}u_{1,k}(u_{k,k})^{-1}F(U_{(k)})\end{equation} hence by Lemma \ref{l11} we get that $$P_1(r,s)=-\sum _{k=2}^ru_{1,k}(u_{k,k})^{-1}P_k(r,s)\ .$$ By the inductive hypothesis we have:
$$P_k(r,p)=\left\{ \begin{array}{cc} 1 & \mbox{\rm if } \ k=r \\ \\ -(u_{r,r})^{-1}F(U_{(k)}^{(r)}) & \mbox{\rm if } \ k\in \{ 2,\ldots ,r-1\} \ \ ; \end{array} \right. $$ hence $$P_1(r,p)=-u_{1,r}(u_{r,r})^{-1}+\sum _{k=2}^{r-1}u_{1,k}(u_{k,k})^{-1} (u_{r,r})^{-1}F(U_{(k)}^{(r)})= $$ $$=-(u_{r,r})^{-1}\bigg[ u_{1,r}-\sum _{k=2}^{r-1}u_{1,k}(u_{k,k})^{-1}F(U_{(k)}^{(r)})\bigg] =-(u_{r,r})^{-1}F(U^{(r)}) \ .$$
Let $s\neq p\,  $. By the inductive hypothesis we have
$$P_k(r,s)=\left\{ \begin{array}{cc}-(u_{s,s})^{-1}F(U_{(s)}) &
\mbox{\rm if }\ k=r \\ \\ (u_{r,r})^{-1}(u_{s,s})^{-1}F(U_{(s)})F(U^{(r)})
& \mbox{\rm if } \ k\in \{ 2,\ldots ,r-1\} \  \
.\end{array}\right. $$ Hence
$$P_1(r,s)=$$ $$=u_{1,r}(u_{r,r})^{-1}(u_{s,s})^{-1}F(U_{(s)})-\sum _{k=2}^{r-1}u_{1,k}(u_{k,k})^{-1}(u_{r,r})^{-1}(u_{s,s})^{-1}F(U_{(s)})F(U_{(k)}^{(r)})=$$ $$= (u_{s,s})^{-1}(u_{r,r})^{-1}F(U_{(s)})\bigg[ u_{1,r}-\sum _{k=2}^{r-1}u_{1,k}(u_{k,k})^{-1}F(U_{(k)}^{(r)})\bigg] = $$ $$= (u_{s,s})^{-1}(u_{r,r})^{-1}F(U_{(s)})F(U_{(r)})\ .\qquad \square \vspace{2mm} $$
\begin{corollary}\label{symmetry} We have that $F(U^{(2)})=u_{1,2}$ and
$$F(U^{(h)})=u_{1,h}-\sum _{k=2}^{h-1}u_{k,h}(u_{k,k})^{-1}F(U^{(k)})$$ for $h=3,\ldots, p$.
\end{corollary}
\pf The first claim is obvious; if $p>2$ by Lemma \ref{coeff} we get that
$$F(U)=u_{1,p}-\sum _{k=2}^{p-1}u_{k,p}(u_{k,k})^{-1}F(U^{(k)})\ ; $$ this result can be generalized to the matrices $U^{(h)}$. \hspace{5mm} $\square $ \vspace{2mm}\newline
For $k=1,\ldots ,p-1$ let $\widehat {U_{(k)}}$ be the submatrix of $U_{(k)}$ obtained by cancelling the first column and the last row. The following Proposition describes $F(U)$ more clearly.
\begin{proposition}\label{determinant} If $U$ is a submatrix of $\Xi ^{\Phi }$ with the property $\star $) then $$F(U)=(-1)^p\, \mbox{\rm det }\widehat U \, \bigg( {\displaystyle \prod _{i=2}^{p-1}u_{i,i}} \bigg) ^{-1} \ \ .$$ \end{proposition}
\pf If $p=2$ the claim is true, since $F(U)=u_{1,2}$. Hence we can prove the claim by induction on $p$. We have that $F(U)$ is expressed by equality \ref{formaF}, moreover by the inductive hypothesis we have that $$F(U_{(k)})=(-1)^{p-k+1}\, \mbox{\rm det }\widehat {U_{(k)}} \, \bigg( {\displaystyle \prod _{i=k+1}^{p-1}u_{i,i}}\bigg) ^{-1} \ \ .$$ Hence we get
$$F(U)=u_{1,p}-\sum _{k=2}^{p-1}u_{1,k}\, (-1)^{p-k+1}\, \mbox{\rm det }\widehat {U_{(k)}} \, \bigg( {\displaystyle \prod _{i=k}^{p-1}u_{i,i}}\bigg) ^{-1} $$
which implies that $$F(U)\, {\displaystyle \prod _{h=2}^{p-1} u_{h,h}}=u_{1,p}\, {\displaystyle \prod _{h=2}^{p-1} u_{h,h}}-\sum _{k=2}^{p-1}u_{1,k}\, (-1)^{p-k+1}\, {\displaystyle \prod _{h=2}^{k-1} u_{h,h}}\, \mbox{\rm det }\widehat {U_{(k)}}\ \ .$$
For $k=2,\ldots ,p-1$ we denote by $D_k$ the determinant of the submatrix of $\widehat U$ obtained by cancelling the row of index $1$ and the column of index $k-1$ (which in $U$ has index $k$); by the previous equality we get $$F(U)\, {\displaystyle \prod _{h=2}^{p-1} u_{h,h}}=u_{1,p}\, {\displaystyle \prod _{h=2}^{p-1} u_{h,h}}+\sum _{k=2}^{p-1}u_{1,k}\, (-1)^{p-k}\, D_k\ .$$ If $p$ is even this is $\mbox{\rm det }\widehat U$ (since $u_{1,k}$ is the entry of indices $(1,k-1)$ in $\widehat U$ for $k=2,\ldots ,p$); if $p$ is odd this is the opposite of $\mbox{\rm det }\widehat U$. \hspace{4mm} $\square $\vspace{2mm}\newline
The previous definition of $F(U)$ can be extended to all the matrices $U$ with property $\star $) (not only submatrices of $\Xi ^{\Phi }$). \newline
Given two $n\times n$ matrices, we will say that a submatrix of one of them corresponds to a submatrix of the other one if they are obtained by choosing the same indices of the rows and of the columns. Let $\Upsilon =(\upsilon _{i,j})$, $\, i,j=1,\ldots ,n $ be a rational function from an affine variety
$\cA $ to $\cN _{\Phi }$. A polynomial in the entries of $\Xi ^{\Phi }$ or a submatrix of $\Xi ^{\Phi }$ can be considered as a map from the set $N_{\Phi }$; hence if $U$ is a square submatrix of $\Xi ^{\Phi }$ the symbol $\, U\, \circ\, \Upsilon \, $ denotes the submatrix of $(\upsilon _{i,j})$ which corresponds to $U$, while $F(U) \, \circ \, \Upsilon $ denotes the rational function on $\cA $ obtained by replacing $\xi _{i,j}$ with $\upsilon _{i,j}$ in $F(U)$.
\newline By induction on the order of $U$ one can easily prove the following result. \begin{proposition}\label{monomials} If $\, \upsilon _{i,j}\, $ is $\, \xi _{i,j}\, $  or $\, 0\, $ for $\, i,j=1,\ldots ,n\, $, $\, U\, $ is a submatrix of $X$ with property $\star $) such that $\, F(U)\, \circ\, \Upsilon =0\, $ and $M$ is one of the fractional monomials whose sum is $F(U)$ then $\, M\, \circ\, \Upsilon =0\, $. \end{proposition}
The definition of $F(U)$ can be generalized to all the upper triangular submatrices of $\Xi ^{\Phi }$; this generalization, which we are going to explain, together with proposition \ref{determinant'} and Corollary \ref{multilinear}, is not used in the following of the paper. Let $U=\Big( u_{i,j}\Big) $, $i,j=1,\ldots ,p$  be an upper triangular submatrix of $\Xi ^{\Phi }$ and let $I_U$ be the subset of $\{ 1, \ldots ,p\} $ of all the elements $l$ such that $u_{l,l}\neq 0$. If $I_U\neq \emptyset $,  let $u\in N$ and $I_l$, $l=1,\ldots ,u$ be such that:  \begin{itemize} \item[1)] $I_l\subseteq \{ 1, \ldots ,p\} $ and the difference between the minimum and the maximum of $I_l$ is the cardinality of $I_l$; \item[2)] the submatrix $U_l$ of $U$ obtained by choosing $I_l$ as set of indices of the rows and the columns has the property $\star $) and $u$ is the maximum possible. \end{itemize}
If $I_U=\emptyset $ we set $F(U)=0$, otherwise we set $F(U)= {\displaystyle \prod _{l=1}^u} F(U_l)$. We define the matrix $\widehat U$ associated to $U$ as before. The matrix $U$ is conjugated to a matrix such that for $i\in I_U$ the elements $u_{i,k},\ u_{k,i}$ are $0$ for $k\in \{ 1,\ldots ,p\} - \{ i\} $, hence by Proposition \ref{determinant} we get the following result.
\begin{proposition} \label{determinant'} If $U$ is an upper triangular matrix then $$F(U)=(-1)^p\, \mbox{\rm det }\widehat U \, \bigg( {\displaystyle \prod _{i\in I_U}u_{i,i}}\bigg) ^{-1} \ \ .$$ \end{proposition}
We observe the following result, which is obvious by the multilinearity of the determinant and the definition of $F(U)$.
\begin{corollary} \label{multilinear} If $l\in \{ 2, \ldots ,p\} $, $U'$ is obtained from $U$ by replacing $u_{h,l}$ with $u'_{h,l} $ for $h=1,\ldots , l$ and $\widetilde U$ is obtained from $U$ by replacing $u_{h,l}$ with $u_{h,l} +u'_{h,l} $ for $h=1,\ldots , l$ then $F\big( \widetilde U\big)$ can be determined by $F(U)$ and $F\big( U'\big) $ in the following ways:
\begin{itemize} \item[i)] if $u_{l,l}=u'_{l,l}= 0$ then $F\big( \widetilde U\big)\, =\, F(U)\, +\, F(U')\, $;
\item[ii)] if $u_{l,l}=0$ and $u'_{l,l}\neq 0 $ then $$F\big( \widetilde U\big)\, = \,  F(U) \, \big( u'_{i,i}\big) ^{-1} \, +\, F(U') \ ;$$
\item[iii)] if $u_{l,l} =- u'_{l,l}\neq 0$ then $F\big( \widetilde U\big)\, = \, F(U)\, u_{l,l}\, +\, F(U')\, u'_{l,l}$;
\item[iv)] if $u_{l,l}, \ u'_{l,l}\, \neq \, 0$ and $u_{l,l}+ u'_{l,l}\neq 0$ then $$F\big( \widetilde U\big)\, =\,
  F(U)\, u_{l,l} \, \big( u_{l,l}+u'_{l,l}\big) ^{-1}\, + \,  F(U')\, u'_{l,l} \, \big( u_{l,l}+u'_{l,l}\big) ^{-1}\ .$$
\end{itemize}
\end{corollary}
\subsection{On a type of subspace of $\cN $ and one of its open subsets}\label{special}
For some maps $\Phi $ the maximum nilpotent orbit which intersects $N_{\Phi }$ contains the whole open subset $\cN _{\Phi }$; this is explained in the following Proposition.
 \begin{proposition}
\label{Phi} If $\; \Phi $ is not decreasing and the restriction of
$\; \Phi $ to the subset $\Phi ^{-1}\big( \{ 2,\ldots ,n\} \big) $ is
increasing there exists a nilpotent orbit which contains $\cN _{\Phi }$.
\end{proposition} \pf For $k\in \N $ we can define the map $\Phi ^k$
from $\{ 1,\ldots ,n\} $ to $\{ 2,\ldots ,n+1\} $ by induction on
$k$, as follows: $\Phi ^1=\Phi $ and $$ \Phi ^k(i)=\left\{
\begin{array}{ll} n+1 & \mbox{\rm if $\Phi ^{k-1}(i)=n+1$} \\ \\
\Phi \big( \Phi ^{k-1}(i)\big) & \mbox{\rm otherwise} \
.\end{array} \right. $$ If $\Phi $ is not decreasing and the
restriction of $\Phi $ to $\Phi ^{-1}\big( \{ 2,\ldots ,n\} \big)
$ is increasing then $\Phi ^k$ is not decreasing and the
restriction of $\Phi ^k$ to $\big( \Phi ^k\big) ^{-1}\big( \{
2,\ldots ,n\} \big) $ is increasing for all $k\in \N $; moreover
the rank of all the elements of $\cN_ {\Phi ^k}$  is the number of
their nonzero rows, that is $$\Big| \{ i\in \{1,\ldots ,n\} \ |\
\Phi ^k(i)\neq n+1\} \Big| \ .$$ If $\Xi ^{\Phi }\in \cN _{\Phi } $ then $\Big( \Xi ^{\Phi }\Big) ^k\in
\cN _{\Phi ^k}$  for all $k\in \N $; this can be proved by induction
on $k$ and considering the transpose matrices of $\Xi ^{\Phi }$ and
$\Big( \Xi ^{\Phi }\Big) ^{k-1}$. Hence we get that rank $\Big( \Xi ^{\Phi }\Big) ^k$ is the same for all $\Xi ^{\Phi }\in
\cN _{\Phi } $, which proves the claim.\hspace{4mm}$\square $
\subsection{A characterization of the maximum nilpotent orbit intersecting a subvariety of $\cN $}\label{characterization}
As in Subsection \ref{special}, let $\Phi $ be a map from $\{ 1,\ldots
,n\} $ to $\{ 2,\ldots ,n+1\} $ such that $\Phi (i)>i$ for $i=1,\ldots ,n$. We will say that an element $i$ of $\{ 1,\ldots ,n-1\} $ is $\Phi $ - regular if  $\Phi (i)<\Phi (i+1)$. Similarly, if $\Psi $ is a map from $\{ 1,\ldots
,n\} $ to $\{ 0,\ldots ,n-1\} $ such that $\Psi (i)<i$ for $i=1,\ldots ,n$ we will say that an element $i$ of $\{ 2,\ldots ,n\} $ is $\Psi $ - regular if  $\Psi (i)>\Psi (i-1)$. Let $\Xi ^{\Phi }=(\xi _{i,j})$, $i,j=1,\ldots ,n$ be the matrix of the coordinates of $N_{\Phi }$. Let $\Xi ^{\Phi ,0}$ be the $n\times (n+1)$ matrix whose first $n$ columns are the columns of $\Xi ^{\Phi }$ and whose last column is $0$. If $U$ is a square submatrix of $\Xi ^{\Phi, 0}$ of order $p$, for $h=1,\ldots ,p$ let $i_U(h)$ be the index of the $h$-th row of $\Xi ^{\Phi ,0}$ chosen for $U$ and let $j_U(h)$ be the index of the $h$-th column of $\Xi ^{\Phi ,0}$ chosen for $U$. If $\, \Phi $ isn't not decreasing or the restriction of
$\, \Phi $ to the subset $\Phi ^{-1}\big( \{ 2,\ldots ,n\} \big) $ is not
increasing, let $U^{\Phi }$ be the square submatrix of $\Xi ^{\Phi ,0}$ of order  $p_{\Phi }>1$ with  property $\star $) and the following properties:
\begin{itemize} \item[ $\star _1$)] $i_{U^{\Phi }}(h+1)=i_{U^{\Phi }}(h)+1$, $\, j_{U^{\Phi }}(h+1)=j_{U^{\Phi }}(h)+1$ for $h=1,\ldots , p_{\Phi}-2$;
 \item[ $\star _2$)] $\ j_{U^{\Phi }}(h)=\Phi (i_{U^{\Phi }}(h))$ for $h=2,\ldots ,p_{\Phi }-1$; \item[ $\star _3$)] $\ \Phi (i_{U^{\Phi }}(1))\geq \Phi (i_{U^{\Phi }}(2))$ (i.e. $i_{U^{\Phi }}(1)$ is not $\Phi $ - regular) and $\ j_{U^{\Phi }}(p_{\Phi })< \Phi \Big( i_{U^{\Phi }}(p_{\Phi })\Big) $ or both are $n+1$ (i.e. if $\Psi $ is such that $\cN ^{\Psi }=\cN _{\Phi }$ then $j_{U^{\Phi }}(p_{\Phi })$ is not $\Psi $ - regular);
\item[$\star _4$)] $j_{U^{\Phi }}(p_{\Phi })$ is the minumum element of $\{ 1,\ldots ,n\} $ for which there is a submatrix $U^{\Phi }$ of $\Xi ^{\Phi ,0}$ with properties $\star _i$), $i=1,\ldots ,4$.
\end{itemize}
Let $G_{\Phi }$ be an $n\times n $ upper triangular matrix over the field of the rational functions of $N_{\Phi }$ such that for $i,j=1,\ldots ,n$ the entry of $G_{\Phi }$ of indices $(i,i)$ is $1$ and the entry of $G_{\Phi } $ of indices $(i,j)$ is $0$ if $(i,j)\neq (i_{U^{\Phi }}(1),i_{U^{\Phi  }}(h))$ for all $h\in \{ 2,\ldots ,p_{\Phi }-1\} $. The entry of $G_{\Phi }$ of indices $(i_{U^{\Phi }}(1), i_{U^{\Phi }}(h))$ will be denoted by $g_{\Phi, h}$ for $h=2,\ldots ,p_{\Phi }-1$. If $j\in \{ j_{U^{\Phi }}(2),\ldots ,n\} $ the entry of $(G_{\Phi })^{-1}\Xi ^{\Phi }G_{\Phi }$ of indices $(i_{U^{\Phi }}(1),j)$ is
$${\displaystyle \xi _{i_{U^{\Phi }}(1),j}-\sum _{k=2}^{p_{\Phi }-1}\xi _{i_{U^{\Phi }}(k),j}\, g_{\Phi, k} } = {\displaystyle \xi _{i_{U^{\Phi }}(1),j}-\sum _{k\in \{ 2,\ldots ,p_{\Phi }-1\} ,\ j_{U^{\Phi }}(k)\leq j}\xi _{i_{U^{\Phi }}(k),j}\, g_{\Phi, k}} \, .$$
For $h=2,\ldots ,p_{\Phi }-1$ we define $g_{\Phi, h}$ by induction on $h$ in the following way: we set $g_{\Phi, 2}=\xi _{i_{U^{\Phi }}(1),j_{U^{\Phi }}(2)}(\xi _{i_{U^{\Phi }}(2),j_{U^{\Phi  }}(2)})^{-1}$ and if $h\in \{ 3,\ldots ,p_{\Phi }-1\} $ we set \begin{equation}\label{eq2} g_{\Phi, h}=\bigg( \xi _{i_{U^{\Phi }}(1),j_{U^{\Phi }}(h)}-\sum _{k=2}^{h-1}\xi _{i_{U^{\Phi }}(k),j_{U^{\Phi }}(h)}\, g_{\Phi, k}\bigg) (\xi _{i_{U^{\Phi }}(h),j_{U^{\Phi }}(h)})^{-1}\ ,\end{equation} hence $g_{\Phi, h}$ is such that the entry of $(G_{\Phi })^{-1}\Xi ^{\Phi }G_{\Phi }$ of indices $(i_{U^{\Phi }}(1), j_{U^{\Phi }}(h))$ is $0$.
For $i,j\in \{ 1,\ldots ,n\} $ let $\, \xi'_{i,j}\, $ be the entry of $\, (G_{\Phi })^{-1}\Xi ^{\Phi }G_{\Phi }\, $ of indices $(i,j)$; by the definition of $G_{\Phi }$ we get the following result.
\begin{proposition}\label{afterdet} The entry $\, \xi'_{i_{U^{\Phi }}(1),j_{U^{\Phi }}(p_{\Phi })}\, $ of $\, (G_{\Phi })^{-1}\Xi ^{\Phi }G_{\Phi }\, $ having indices $\, \Big( i_{U^{\Phi }}(1),j_{U^{\Phi }}(p_{\Phi })\Big) \,$ is the rational function $F\, \Big( U^{\Phi }\Big) \, $.\end{proposition}
\pf By Proposition \ref{determinant} we have that $$F\Big( U^{\Phi } \Big) =(-1)^{p_{\Phi }}\, \mbox{\rm det }\widehat {U^{\Phi } }\, \bigg( {\displaystyle \prod _{h=2}^{p_{\Phi }-1}\, \xi_{i_{U^{\Phi }}(h),j_{U^{\Phi }}(h)}}\bigg) ^{-1 } \ \ .$$ Since the determinant of $\, \widehat {U^{\Phi }}\, $ is the same as the determinant of the submatrix of $\, (G_{\Phi })^{-1}\Xi ^{\Phi }G_{\Phi }\, $ which corresponds to $\, \widehat {U^{\Phi }}\, $ we get that:
$$\mbox{\rm det }\widehat {U^{\Phi }}= (-1)^{p_{\Phi }}\, \xi'_{i_{U^{\Phi }}(1),j_{U^{\Phi }}(p_{\Phi })}\, {\displaystyle \prod _{h=2}^{p_{\Phi }-1}\, \xi_{i_{U^{\Phi }}(h),j_{U^{\Phi }}(h)}}\ ,$$ which implies the claim. \hspace{4mm} $\square $
\begin{example}\label{ex2}{\rm If $U^{\Phi }$ is the matrix $U$ of Example \ref{ex1} where $u_{1,1}=u_{5,5}=0$ then $\xi _{i,j}=u_{i,j}$ for $i=1,\ldots ,4$ and $g_{\Phi ,2}=u_{1,2}(u_{2,2})^{-1}=F(U^2)(u_{2,2})^{-1}$, $$g_{\Phi, 3}=(u_{1,3}-u_{2,3}g_{\Phi , 2})(u_{3,3})^{-1}=F(U^3)(u_{3,3})^{-1}\ ,$$ $$g_{\Phi ,4}=(u_{1,4}-u_{2,4}g_{\Phi , 2}-u_{3,4}g_{\Phi ,3})(u_{4,4})^{-1}=F(U^4)(u_{4,4})^{-1}\ ;$$
hence $$\xi '_{1,j}=\xi _{1,j}-\xi _{2,j}g_{\Phi , 2}-\xi _{3,j}g_{\Phi ,3}-\xi _{4,j}g_{\Phi ,4}=$$ $$=\xi _{1,j}-\xi _{2,j}(u_{2,2})^{-1}F(U^2)-\xi _{3,j}(u_{3,3})^{-1}F(U^3)-
\xi _{4,j}(u_{4,4})^{-1}F(U^4)=F(U)\ .$$}\end{example}
If $j_{U^{\Phi }}(p_{\Phi})< n$, for $j=j_{U^{\Phi }}(p_{\Phi})+1, \ldots ,n$ let $\Big( U^{\Phi }\Big) _j$ be the submatrix of $\Xi ^{\Phi }$ obtained by choosing the rows of $\Xi ^{\Phi }$ of indices $i_{U^{\Phi }}(1), \ldots , i_{U^{\Phi }}(p_{\Phi })$ and the columns of $\Xi ^{\Phi }$ of indices $j_{U^{\Phi }}(1), \ldots , j_{U^{\Phi }}(p_{\Phi }-1),j$.\newline By the definition of $G_{\Phi }$ and by Proposition \ref{afterdet} we get the following result.
\begin{corollary}\label{F(U)} The matrix $(G_{\Phi })^{-1}\Xi ^{\Phi }G_{\Phi }$ has the following entries: \begin{itemize} \item[1)] the entry of indices $(i_{U^{\Phi }}(1),j_{U^{\Phi }}(k))$ is $0$ for $k=2, \ldots , p_{\Phi } -1$, the entry of indices $(i_{U^{\Phi }}(1),j_{U^{\Phi }}(p_{\Phi }))$ is $F(U^{\Phi })$, if $j_{U^{\Phi }}(p_{\Phi})< n$ the entry of indices $(i_{U^{\Phi }}(1),j)$ is $F\Big( \Big( U^{\Phi }\Big) _j\Big) $ for $j=j_{U^{\Phi }}(p_{\Phi})+1, \ldots ,n$;
\item[2)] the entry of indices $(l, i_{U^{\Phi }}(h))$ is $$\xi _{l,i_{U^{\Phi }}(h)}+F\Big( (U^{\Phi })^{(h)}\Big) (\xi _{i_{U^{\Phi }}(h),j_{U^{\Phi }}(h)})^{-1}\, \xi _{l,i_{U^{\Phi }}(1)}$$ for $h=2,\ldots ,p_{\Phi }-1$ and $l=1,\ldots ,n$;
\item[3)] each one of the other entries is the same as the corresponding entry in $\Xi ^{\Phi }$. \end{itemize}
\end{corollary}
Let us consider the following permutation of $\{ 1,\ldots ,n\} $: $$\sigma _{\Phi }= (i_{U^{\Phi }}(1)\, ,\, i_{U^{\Phi }}(p_{\Phi }-1)\, ,\, i_{U^{\Phi }}(p_{\Phi }-2)\, ,\,  \ldots \, ,\,  i_{U^{\Phi }}(2)) \ ; $$ it induces a permutation of the canonical basis $\{ e_1,\ldots ,e_n\} $ of $K^n$, which is associated to the endomorphism $S_{\Phi }$ of $K^n$ such that $S_{\Phi }(e_i)=e_{\sigma _{\Phi }(i)}$ for $i=1,\ldots ,n$.
If $\, \Phi $ isn't not decreasing or the restriction of
$\, \Phi $ to the subset $\Phi ^{-1}\big( \{ 2,\ldots ,n\} \big) $ is not
increasing,
let $\Sigma  _{\Phi }\ : \ \cN \longrightarrow  \cN $ be defined by $\Xi \ \longmapsto \ \Big( S_{\Phi }\Big) ^{-1}\, \Big( G_{\Phi }\Big) ^{-1}\, \Xi \, G_{\Phi }\, S_{\Phi }$; otherwise let $\Sigma _{\Phi }$ be the identity in $\cN $. \newline
Let $\Upsilon =(\upsilon _{i,j})$, $\, i,j=1,\ldots ,n$ be a rational function from an affine variety
$\cA $ to $\cN $. For any such rational function, let $\Phi _{\Upsilon }$ be the map from $\{ 1,\ldots
,n\} $ to $\{ 2,\ldots ,n+1\} $ such that $\upsilon _{i,j}=0$ if $j\in \{ 1, \Phi _{\Upsilon } (i)-1\} $ and either $\, \upsilon _{i, \Phi _{\Upsilon }(i)}\neq 0$ or $\Phi _{\Upsilon }(i)=n+1$, for all $i=1,\ldots ,n$. Similarly, we can define $\Psi ^{\Upsilon }$ as the map from $\{ 1,\ldots
,n\} $ to $\{ 0,\ldots ,n-1\} $ such that $\upsilon _{i,j}=0$ if $i\in \{\Psi ^{\Upsilon } (j)+1, \ldots ,n\} $ and either $\, \upsilon _{\Psi ^{\Upsilon }(j),j}\neq 0$ or $\Psi ^{\Upsilon }(j)=0$, for all $j=1,\ldots ,n$ (the results of this section will concern the map $\Phi _{\Upsilon }$ rather than the map $\Psi ^{\Upsilon }$).\newline
We consider a certain rational function $\Upsilon $; for all $m\in \N \, \cup \, \{0\} $ we define a rational function $\Upsilon _m$ from $\cA $ to $\cN $ in the following way: $\ \Upsilon_0=\Upsilon $ and $$\ \Upsilon _ m=\Sigma _{\Phi _{\Upsilon _{m-1}}}\, \circ \, \Upsilon _{ m-1} $$ for all $m\in \N $ such that $m\geq 1$.
By these definitions we get the following Lemma and we can easily prove the following Proposition.
\begin{lemma}\label{obv} For all $m\in \N $ the maximum nilpotent orbit which intersects $\Upsilon (\cA )$ is the maximum nilpotent orbit which intersects $\Upsilon _ m (\cA )$. \end{lemma}
\begin{proposition}\label{steps} For any rational function $\Upsilon \, : \, \cA \longrightarrow \cN $ there exists $m\in \N\, \cup \, \{ 0\} $ such that $\Phi_{\Upsilon _{m}}$ is not decreasing and the restriction of $\ \Phi _{\Upsilon _{m}}\ $  to $\ \Big( \Phi_{\Upsilon _{m}}\Big) ^{-1}(\{ 1,\ldots ,n\} )$ is increasing. \end{proposition}
\pf Let $\nu _{m}$ be the cardinality of the set of all $i\in \{ 1,\ldots ,n\} $ such that $\{ l\in \{ i+1,\ldots ,n\}\ | \ \Phi _{\Upsilon _{m}} (l)\neq n+1, \ \Phi _{\Upsilon _{m}} (l)\leq \Phi _{\Upsilon _{m}} (i)\} \neq \emptyset $. The claim is true if there exists $m\in \N \, \cup \, \{ 0\} $ such that $\nu _{m}=0$, hence it is enough to prove that there exists $m\in \N \, \cup \, \{ 0\} $ such that $\nu _m<\nu _0$. If $n-i_{U^{\Phi _{\Upsilon } }}(1)=2$ then $\nu _{ 1} <\nu _{0}$, hence we can prove the claim by induction on $n-i_{U^{ \Phi _{\Upsilon }}}(1)$. By the definition of $\Upsilon _1$ we get that $\nu _{ 1} <\nu _0$ or $\nu _1=\nu _0$ and $n-i_{U^{\Phi _{\Upsilon } }}(1)\, >\, n-i_{U^{\Phi _{\Upsilon , 1 }}}(1)$; hence the claim follows by the inductive hypothesis. \hspace{4mm} $\square $\vspace{2mm} \newline
We will denote by $m_{\Upsilon }$ the minimum of all $m\in \N\, \cup \, \{ 0\} $ which have the property expressed in Proposition \ref{steps}.\begin{example}\label{ex3} {\rm Let $  \Upsilon $ be the upper triangular matrix of order $13$ with the following form:
{ {$$\pmatrix{ & \ast & \ast  & \ast  & \ast  &  \ast  & \ast  & \ast  &
\ast  &  \ast  & \ast  & \ast  &  \ast  \cr &  &\ast  & \ast  & \ast &  \ast  &
\ast  & \ast  & \ast  &  \ast  & \ast  & \ast  &  \ast  \cr & &  & \ast  & \ast
& \ast  & \ast  & \ast  &\ast  & \ast  & \ast  & \ast  &  \ast  \cr & & &  & \ast  &
 { 0} & \ast  & \ast  & \ast &  \ast   & \ast   & \ast  &
\ast  \cr  & & & &  &  { 0} & { 0} & { 0} & \diamond &
 { 0} &  \diamond  &  \diamond  &  \ast  \cr & & & &  &  &  \diamond & \diamond &  \diamond &
  \diamond &  \diamond  &  \diamond  &  \ast  \cr & & & & &  &  &  \diamond &  \diamond &   \diamond  &  \diamond
&  \diamond  &  \ast \cr &  & & & & & & &  \diamond  &   \diamond &  \diamond &  \diamond &
 \ast  \cr &  & & & &   & & &  &  { 0} & \circ  & \circ  &
\circ   \cr &
 & & & &   & & & &   & \circ  & \circ  & \circ  \cr &
& & & &  & & &  &   & & \circ  &  \circ  \cr &  & & & &  & &
& &  & &  &  \circ   \cr &  & & & & & & &  & & & &  \cr  } \
.$$}}
Then $\ i_{U^{\Phi _{\Upsilon }}}(1)=5$, $\ j_{U^{\Phi _{\Upsilon }}}(1)=6$, $\ p_{\Phi _{\Upsilon }}=5$ and:
{ {$$\Upsilon _1=\pmatrix{ & \ast  & \ast  &
\ast  & \ast   & \ast  & \dot{\ast } & \ast  &{\ast }  & {\ast } & {\ast }  & \ast  & \ast \cr &
& \ast  & \ast  & \ast  & \ast & \dot{\ast  } & \ast  & {\ast}  & {\ast } & {\ast }  & \ast  & \ast \cr
& &  & \ast & \ast   & \ast  & \dot{\ast } & \ast   & {\ast } & {\ast } & {\ast }  & \ast  &\ast  \cr
& & &  & { 0} & {\ast } & \dot{\ast } &  \ast &\ast & {\ast } & {\ast } &
{\ast } & \ast   \cr   & & & &   & \diamond
& {\diamond } &  & {\diamond } & {\diamond } & {\diamond }& \diamond  & \ast  \cr & & & &   & & {\diamond } & & {\diamond } &
{\diamond }  & {\diamond }  & \diamond  & \ast  \cr & & & & & & &   & {\diamond } &{\diamond } & {\diamond }
 & \diamond &
 \ast  \cr   & & & & &  &  &  &
& \triangle  & \ddot{\diamond } & \ddot{\diamond } & \ddot{\ast } \cr&  & & & &   & & &  &  { 0} & \circ  & \circ  &
\circ   \cr &
 & & & &   & & & &   & \circ  & \circ  & \circ  \cr &
& & & &  & & &  &   & & \circ  &  \circ  \cr &  & & & &  & &
& &  & &  &  \circ   \cr &  & & & & & & &  & & & &  \cr } \vspace{2mm}$$ }} where the empty spaces correspond to entries equal to $0$. The entry $\triangle $ is $F\Big( U^{\Phi _{\Upsilon  }}\Big) \, \circ \, \Upsilon $; the two dots $\ddot{\ } $ above certain entries denote that it is obtained from the corresponding entry of $\Upsilon $ by adding a rational function according to 1) of Corollary \ref{F(U)}, while the dot $\dot{\ } $ above certain entries denotes that it is obtained from the corresponding entry of $\Upsilon $ by adding to it a rational function according to 2) of Corollary \ref{F(U)}.  Now we have $\ i_{U^{\Phi _{\Upsilon _1}}}(1)=4$, $\ j_{U^{\Phi _{\Upsilon _1 }}}(1)=5$, $\ p_{\Phi _{\Upsilon _1}}=4$.  If $\Upsilon $ is a generic rational function then:\vspace{2mm}
{ {$$\Upsilon _2=\pmatrix{ & \ast  & \ast  &
\dot{\ast } & \dot{\ast }  & \ast  & \dot{\ast } & \ast  &{\ast }  & {\ast } & {\ast }  & \ast  & \ast \cr &
& \ast  & \dot{\ast } & \dot{\ast } & \ast & \dot{\ast  } & \ast  & {\ast}  & {\ast } & {\ast }  & \ast  & \ast \cr
& &  & \dot{\ast } & \dot{\ast }  & \ast  & \dot{\ast } & \ast   & {\ast } & {\ast } & {\ast }  & \ast  &\ast  \cr
   & & &   & {\diamond }
&  & {\diamond } & & {\diamond } & {\diamond } & {\diamond }& \diamond  & \ast  \cr & & &    & & & {\diamond } & & {\diamond } &
{\diamond }  & {\diamond }  & \diamond  & \ast  \cr & &   &  &  & &  & \nabla & \ddot{\ast } &\ddot{\ast } & \ddot{\ast } & \ddot{\ast } &
\ddot{\ast }   \cr & & & & & & &   & {\diamond } &{\diamond } & {\diamond }
 & \diamond &
 \ast  \cr   & & & & &  &  &  &
& \triangle  & \ddot{\diamond } & \ddot{\diamond } & \ddot{\ast } \cr &  & & & &   & & &  &  { 0} & \circ  & \circ  &
\circ   \cr &
 & & & &   & & & &   & \circ  & \circ  & \circ  \cr &
& & & &  & & &  &   & & \circ  &  \circ  \cr &  & & & &  & &
& &  & &  &  \circ   \cr &  & & & & & & &  & & & &  \cr }\ . \vspace{2mm}$$ }}
Since $\ i_{U^{\Phi _{\Upsilon _2}}}(1)=9$, $\ j_{U^{\Phi _{\Upsilon _2 }}}(1)=10$ and $\ p_{\Phi _{\Upsilon _2}}=5$, we have
{ {$$\Upsilon _3=\pmatrix{ & \ast  & \ast  &
\dot{\ast } & \dot{\ast }  & \ast  & \dot{\ast } & \ast  &{\ast }  & \dot{\ast } & \dot{\ast }  & \dot{\ast } & \ast \cr &
& \ast  & \dot{\ast } & \dot{\ast } & \ast & \dot{\ast  } & \ast  & {\ast}  & \dot{\ast } & \dot{\ast }  & \dot{\ast } & \ast \cr
& &  & \dot{\ast } & \dot{\ast }  & \ast  & \dot{\ast } & \ast   & {\ast } & \dot{\ast } & \dot{\ast }  & \dot{\ast } &\ast  \cr
   & & &   & {\diamond }
&  & {\diamond } & & {\diamond } & \dot{\diamond } & \dot{\diamond }& \dot{\diamond } & \ast  \cr & & &    & & & {\diamond } & & {\diamond } &
\dot{\diamond }  & \dot{\diamond }  & \dot{\diamond } & \ast  \cr & &   &  &  & &  & \nabla & \ddot{\ast } &\dot{\ddot{\ast }}& \dot{\ddot{\ast }} & \dot{\ddot{\ast } } &
\ddot{\ast }   \cr & & & & & & &   & {\diamond } & \dot{\diamond } & \dot{\diamond }
 & \dot{\diamond } &
 \ast  \cr   & & & & &  &  &  &
& {\triangle } & \ddot{\diamond } & \ddot{\diamond } & \ddot{\ast } \cr &
 & & & &   & & & &   & {\circ } & {\circ } & \circ  \cr &
& & & &  & & &  &   & & {\circ } &  \circ  \cr &  & & & &  & &
& &  & &  &  \circ   \cr &  & & & & & & &  & & & &  \cr &  & & & &   & & &  &   &   &   &
   \cr }\ . \vspace{2mm}$$ }}
Hence in this case, if $\Upsilon $ is a generic rational function, $m_{\Upsilon }$ is 3.}\end{example}
By Proposition \ref{determinant} and property $\star _4$) we get the following result.
\begin{proposition}\label{equalphi} If $\Upsilon $ and $\Upsilon '$ are two rational functions from the same affine variety
$\cA $ to $\cN $ such that any submatrix of $\Upsilon $ has the same rank as the corresponding submatrix of $\Upsilon '$ then $m_{\Upsilon }=m_{\Upsilon '}$ and $\Phi _{\Upsilon _{m_{\Upsilon }}}=\Phi _{\Upsilon ' _{m_{\Upsilon '}}}$. \end{proposition}
\pf For all $m\in \N $ we have that $F\Big( U^{\Phi _{\Upsilon _m }}\Big) \, \circ \, \Upsilon _m$ is $0$ iff $\widehat{U^{\Phi _{\Upsilon _m }}}\, \circ \, \Upsilon _m$ is singular, and this matrix is obtained from a submatrix $U$ of $\Upsilon $ by the operations described in Corollary \ref{F(U)}. If $\widehat{U^{\Phi _{\Upsilon _m }}}\, \circ \, \Upsilon _m$ has a row which is obtained by row operations according to 2) of Corollary \ref{F(U)} then the rows which are involved in these operations are still rows of $\widehat{U^{\Phi _{\Upsilon _m }}}\, \circ \, \Upsilon _m$, since they are some of the $\Phi  _{\Upsilon _m }$ - regular rows just above that row. If instead $\widehat{U^{\Phi _{\Upsilon _m }}}\, \circ \, \Upsilon _m$ has a column which is obtained by column operations according to 1) of Corollary \ref{F(U)} then by  $\star _4$) there are two possibilities for a column which is involved in these operations: it is the last column of $\widehat{U^{\Phi _{\Upsilon _m }}}\, \circ \, \Upsilon _m$ or it is the last column of $\widehat{U^{\Phi _{\Upsilon _{m'} }}}\, \circ \, \Upsilon _{m'}$ for $m'<m$; in this last case $\ F\Big( U^{\Phi _{\Upsilon _{m'}}}\Big) \, \circ \, \Upsilon _{m'}\ \neq \ 0\ $, which implies that $F\Big( U^{\Phi _{\Upsilon _m }}\Big) \, \circ \, \Upsilon _m\ \neq \ 0 \ $. This shows that $\widehat{U^{\Phi _{\Upsilon _m }}}\, \circ \, \Upsilon _m$  is singular iff $U$ is singular.
 \hspace{4mm} $\square $ \vspace{2mm}\newline
By Proposition \ref{Phi} we get the following result.
\begin{corollary}\label{aftersteps} The maximum nilpotent orbit which intersects $\Upsilon (\cA )$ is determined by the map $\Phi _{\Upsilon _{m_{\Upsilon }}}$.\end{corollary}
{\bf Proof of Theorem \ref{equalorbit}} \ The claim is a consequence of Corollary \ref{aftersteps} and Proposition \ref{equalphi}. \hspace{4mm} $\square $
 \subsection{Rational functions into $M(n,K)$ with property $\ast $}\label{P}
Let $\cA $ be an affine variety over $K$ and let $\Upsilon =(\upsilon _{i,j})$, $\, i,j=1,\ldots ,n$ be a rational function from $\cA $ to $M(n,K)$.  We define a $D$ - product of entries of a matrix as a product of entries such that there are not two of them in the same row or column and such that the number of those entries is the order of the matrix. We will say that the rational function $\Upsilon $ has property $\ast $) if for any nonzero submatrix $V$ of $(\upsilon
_{i,j})$ there exist a coordinate $x$ of $\cA $ and only one entry $f$ of $V$ different from $0$ such that:
\begin{itemize} \item[$\ast $)] if the sum of all the $D$ - products of entries of $V$ among which there is $f$ depends on $x$ then the determinant of $V$ also depends on $x$. \end{itemize}
An entry $f$ of $V$ with the property expressed in $\ast $) will be said a $\ast $) - entry of $V$; similarly, a coordinate $x$ of $\cA $  with the property expressed in $\ast $) will be said a $\ast $) - coordinate of $V$. As previously, let $ \Xi =(\xi _{i,j})$ be the generic element of $M(n,k)$ and let $\Upsilon ' =(\upsilon '_{i,j})$ be defined as follows:
$$\upsilon '_{i,j}=\left\{ \begin{array}{lr} \xi _{i,j} & \mbox{\rm if } \ \upsilon _{i,j}\neq 0 \\ 0 & \mbox{\rm otherwise.} \end{array} \right. $$
 If $V$ is a submatrix of $(\upsilon
_{i,j})$ we will denote by $V'$ the submatrix of $(\upsilon '_{i,j})$ which corresponds to $V$.
\begin{proposition}\label{rank} If $\Upsilon $ has
property $\ast $) and $V$ is a submatrix of $\Upsilon =(\upsilon _{i,j})$, $i,j=1,\ldots ,n$ then $\mbox{\rm rank }V\, =\, \mbox{\rm rank } V'$.  \end{proposition}
\pf We can assume that $V$ is a square submatrix of $\Upsilon $; let $q$ be the order of $V$. By property $\ast $) the claim is obviously true if $q$ is $1$, hence we can prove it by induction on $q$. Let us assume that $\mbox{\rm rank }  V' \, =\,  q $, $\, \mbox{\rm rank } V\, <\, q $ and let us prove that this assumption leads to a contradiction.  For $h=1,\ldots ,q$ let $i_V(h)$ be the index of the $h$-th row of $(\upsilon _{i,j})$ chosen for $V$ and let $j_V(h)$ be the index of the $h$-th column of $(\upsilon _{i,j})$ chosen for $V$. We consider $V$ as the matrix of an endomorphism from $\langle e_{j_V(1)},\ldots ,e_{j_V(q)}\rangle $ to $\langle e_{i_V(1)},\ldots ,e_{i_V(q)}\rangle $ over the field $K(\cA )$. Let $h\in \{ 1,\ldots ,q\} $, $k\in \{ 1,\ldots ,q\} $ be such that the entry of $V$ of indices $(h,k)$ is a $\ast $) - entry of $V$. Let $V_{\{ h\} }$ ($\, (V')_{\{ h\} }$) be the submatrix of $V$ (of $V'$) obtained by cancelling the row of $V$ (of $V'$) whose index is $h$, let $V^{\{ k\} }$ ($\, (V')^{\{ k\} }$) be the submatrix of $V$ (of $V'$) obtained by cancelling the column of $V$ (of $V'$) whose index is $k$ and let $V_{\{ h\} }^{\{ k\} }$ be the submatrix of $V_{\{ h\} }$ (of $(V')_{\{ h\} }$) obtained by cancelling the column of $V_{\{ h\} }$ (of $(V')_{\{ h\} }$) whose index is $k$.  If the rank of $(V)_{\{ h\} }^{\{ k\} }$ is $q -1$ then the rank of $V$ is $q$ (since in the determinant of $V$ the coefficient of a $\ast $) - coordinate associated to the $\ast $) - entry of $V$ of indices $(h,k)$ is not $0$). This reduces the proof to the case in which $\mbox{\rm rank }(V)_{\{ h\} }^{\{ k\} }\, < \, q -1$.  By the inductive hypothesis we have that $\mbox{\rm rank } (V')_{\{ h\} }^{\{ k\} }\, < \, q -1$ and $\mbox{\rm rank } (V')_{\{ h\} },\ \mbox{\rm rank }(V')^{\{ k\} } = q -1$. Let us consider $\{ 1,\ldots ,q\} - \{ k\} $ as the set of the indices of $V_{\{ h\} }^{\{ k\} }$ ($(V')_{\{ h\} }^{\{ k\} }$) and let $k'\in \{ 1,\ldots ,q\} - \{ k\} $ be such that the column of $(V')_{\{ h\} }^{\{ k\} }$ whose index is $k'$ is a linear combination over $K(\cA ')$ of the other columns of $(V')_{\{ h\} }^{\{ k\} }$. If we substitute the coefficients and the vectors of this linear combination with their images by $\eta $ we express the column of $V_{\{ h\} }^{\{ k\} }$ whose index is $k'$ as a linear combination over $K(\cA )$ of the other columns of $V_{\{ h\} }^{\{ k\} }$; this follows by the inductive hypothesis, since by the Cramer method the coefficients of a linear combination of vectors which is equal to the zero vector are obtained as quotients of determinants.
We consider a new basis $\{ e'_{j_V(l)},\  l\in \{ 1,\ldots ,q\} - \{ k\} \} $ of $\langle e_{j_V(l)},\ l\in \{ 1,\ldots ,q\} - \{ k\} \rangle $ and a new basis $\{ e'_{i_V(l)},\  l\in \{ 1,\ldots ,q\} - \{ h\} \} $ of $\langle e_{i_V(l)},\ l\in \{ 1,\ldots ,q\} - \{ h\} \rangle $ such that the representation of $(V')_{\{ h\} }^{\{ k\} }$ with respect to these new bases has the column of index $k'$ equal to $0$. We set $e'_{j_V(k)}=e_{j_V(k)}$, $e'_{i_V(h)}=e_{i_V(h)}$ and we denote by $\tilde V$ the representation of $V$ with respect to the bases $\{ e'_{j_V(l)},\  l\in \{ 1,\ldots ,q\} \} $, $\, \{ e'_{i_V(l)},\  l\in \{ 1,\ldots ,q\} \} $. By the previous observation on the rank of $V_{\{ h\} }$ we get that the entry of $\tilde V$ of indices $(h,k')$ is not $0$; by the previous observation on the rank of $V^{\{ k\} }$ we get that the submatrix of $\tilde V$ obtained by cancelling the row of index $h$ and the column of index $k'$ has rank greater than or equal to $q -1$.  But this implies that the rank of $\tilde V$ is $ q $, which is in contradiction with the assumption that the determinant of $V$ is $0$.\hspace{4mm} $\square $
\section{On the nilpotent subalgebras $\cSE _B$ and $ \cSN _B$}\label{centralizer}
\subsection{Description of the varieties ${\cE }_B$, ${\cN }_B$,
$\cSE _B$ and $\cSN _B$} We will consider any $n\times n$ matrix
$X$ as a block matrix $(X_{h,k})$, where $X_{h,k}$ is a $\mu
_h\times \mu _k$ matrix and $h,k=1,\ldots ,t$. Let $ \cU _B$ be
the subalgebra of $M(n,K)$ of all $X$ such that for $1\leq k\leq
h\leq t$ the blocks $X_{h,k}$ and $X_{k,h}$ are upper triangular,
that is have the following form:\vspace{2mm}
$$X_{h,k}=\pmatrix{0 &\ldots & 0 & x_{h,k}^{1,1} & x_{h,k}^{2,1} & \ldots &
x_{h,k}^{\mu _h,1} \cr \vdots &  &  & 0 & x_{h,k}^{1,2} & \ddots &
x_{h,k}^{\mu _h-1,2} \cr \vdots &  &  &  & \ddots & \ddots &
\vdots \cr 0 & \ldots & \ldots & \ldots & \ldots & 0 &
x_{h,k}^{1,\mu _h} \cr},\vspace{2mm}$$
$$X_{k,h}=\pmatrix{x_{k,h}^{1,1} & x_{k,h}^{2,1} & \ldots & x_{k,h}^{   \mu _h,1} \cr
0 & x_{k,h}^{1,2} & \ddots & x_{k,h}^{\mu _h-1,2} \cr \vdots &
\ddots & \ddots & \vdots \cr \vdots &  & 0 & x_{k,h}^{1,\mu _h}
\cr \vdots & & & 0 \cr \vdots &  &  & \vdots \cr 0 & \ldots &
\ldots & 0 \cr}\vspace{2mm}$$ where for $\mu _h=\mu _k$ we omit
the first $\mu _k-\mu _h$ columns and the last $\mu _k-\mu _h$
rows respectively. \newline For $X\in \cU _B $, $\ i,j \in \{
1,\ldots ,u\} $ and $l\in \{ 1,\ldots ,\mbox{\rm min }\{ \mu
_{q_{i }}, \mu _{q_{j }}\} \} $ let
$$X (i,j ,l)=(x_{h,k}^{1,l})\ ,\qquad  q_{i
-1}+1\leq h\leq  q_{i},\ q_{j -1}+1\leq k\leq  q_{j}.$$ We set
$X(i,l)=X(i,i,l)$. Lemma \ref{algebra} for the algebra $\cU _B$
becomes more precise, as follows.
\begin{lemma} \label{1}
For $X\in \cU _B $ we have that: \begin{itemize}  \item[a)] there
exists $G\in $ GL $(n,K)$ such that $G^{-1}X G\,\in \,\cU _B$ and $ (G^{-1}X G)\; (i ,l) $ is
lower (upper) triangular  for $\, i =1,\ldots ,u\, $ and $\, l=
1,\ldots ,\mu _{q_{i }}\, $; \item[b)]
 $X$ is nilpotent iff $\, X( i , l)\, $ is nilpotent for $\, i =1,\ldots ,u\, $ and $\, l= 1,\ldots
,\mu _{q_{i }} \, $.\end{itemize}
\end{lemma}
\pf  We can construct a semisimple subalgebra of $\cU _B $ whose
direct sum with the Jacobson radical of $\cU _B$ is $\cU _B $; the
construction is as follows. For $l=1,\ldots ,\mu _{q_1}$ let
$$U^l=\langle v_{\mu _{q_i},j}^l \ : \ i=1,\ldots ,u,\
j=q_i-q_{i-1},\ldots ,1,\  \mu _{q_i}\geq l\rangle ;$$ then
$K^n={\displaystyle \bigoplus _{l=1}^{\mu _{q_1}}} U^l$ and
$X(U^l)\subseteq {\displaystyle \bigoplus _{i=l}^{\mu
_{q_1}}}U^i$. For $v\in K^n$ let $v={\displaystyle \sum
_{l=1}^{\mu _{q_1}}}v^{(l)}$ where $v^{(l)}\in U^l$ and let
$L_{X,l}\, :\, U^l\to U^l$ be defined by $L_{X,l}(v)=X(v)^{(l)}$.
Then $X$ is nilpotent iff $L_{X,l}$ is nilpotent for $l=1,\ldots
,\mu _{q_1}$. For $l=1,\ldots ,\mu _{q_1}$ let $i _l\in \{
1,\ldots ,u \} $ be such that $l\leq \mu _{q_{i_l}}$ and $\mu
_{q_{i_l+1}}<l$ if $i_l\neq u$. Then the matrix of $L_{X,l}$ with
respect to the basis $\{ v_{\mu _{q_i},j}^l\ \  : \ i=1,\ldots
,u,\ \mu _{q_i}\geq l\} $ is the lower triangular block matrix
$\big( X (i,j ,l)\big) $, $i ,j =1,\ldots ,i _l$, which is
nilpotent iff $X (i, l)$ is nilpotent for $i =1,\ldots ,i _l$. For
$v\in K^n$ let $v={\displaystyle \sum _{i = 1}^{u } v_{(i)}}$
where $v_{(i )}\in \langle v_{\mu_{q_i},j}^l \ :\
j=q_i-q_{i-1},\ldots ,1,\   l=1,\ldots ,u_{q_{i}}\rangle $. For $i
=1,\ldots ,u $ and $l=1,\ldots ,u_{q_{i }}$ let $U^l_{i}=\langle
v_{\mu_{q_i},j}^l\ :\ j=q_i-q_{i-1},\ldots ,1\rangle $ and let $
L_{X,i ,l}\ :\ U^l_{i}\to U^l_{i} $ be defined by $L_{X,i
,l}(v)=L_{X,l}(v)_{(i)}$. Then $ X(i ,l )$ is the matrix of $L
_{X,i,l}$ with respect to the basis $\{ v_{\mu_{q_i},j}^l\ :\
j=q_i-q_{i-1},\ldots ,1\} $ . We can substitute this basis with
another basis of the same subspace such that $X(i,l)$ is upper
triangular, for $i=1,\ldots ,u$ and $j=q_i-q_{i-1},\ldots
,1$.\hspace{4mm} $\square$ \vspace{2mm}\newline
 We
will denote by $ \cSU   _B $ the subspace of $\cU _B $ of all $X$
such that $X(i ,l)$ is lower triangular for $i =1,\ldots ,u$ and
$l=1,\ldots ,\mu _{q_{i }}$. Moreover we will denote by $\cE _{B
}$ the subset of $\cU _B $ of all the nilpotent matrices and we
will set
$$ \cSE _{B }=\cSU   _B\; \cap \; \cE _{B }\ .$$
\begin{lemma} \label{2}
 The centralizer $\cC _B$ of $J$ has the following properties: \begin{itemize} \item[i)] it is the
subspace of  $\cU _B$ of all $X$ such that
$$x_{h,k}^{l,1}=x_{h,k}^{l,2}=\cdots =x_{h,k}^{l,\mu _h+1-l}\ ,
\quad
 \ x_{k,h}^{l,1}=x_{k,h}^{l,2}=\cdots = x_{k,h}^{l,\mu _h+1-l}$$
for $\ 1\leq k\leq h\leq t\ $ and $\ l=1,\ldots ,\mu _h$;
\item[ii)] if $X\in \cC _B$ we can choose $G$ with the property
expressed in a) of lemma \ref{1} and such that
$GJ=JG$.\end{itemize}
\end{lemma}
\pf For i) see \cite{Ait} or Lemma 3.2 of \cite{Bas'}. Using the
notations of the proof of lemma \ref{1}, for $i=1,\ldots ,u$ let
${\big (} c^{(i)}_{h,k}{\big )} $, $h,k=q_i-q_{i-1},\ldots ,1$ be
a $q_i-q_{i-1}$ matrix over $K$ such that the vectors
$$w_{\mu _{q_i},j}^1={\displaystyle \sum
_{k=q_i-q_{i-1}}^1c^{(i)}_{j,k}v_{\mu _{q_i},k}^1}$$  form a basis
with respect to which $L_{X,i,1}$ is upper triangular. If we set
$$w_{\mu _{q_i},j}^l={\displaystyle \sum
_{k=q_i-q_{i-1}}^1c^{(i)}_{j,k}v_{\mu _{q_i},k}^l}$$ for
$l=1,\ldots \mu _{q_i}$ and for $i=1,\ldots ,u$ we get the basis
required by ii).\hspace{4mm} $\square $\vspace{2mm}\newline We can
shortly say that
 $X\in \cC _B$ iff its blocks are upper triangular
Toeplitz matrices.
 By lemma \ref{2} if $\, A\in \cC _B\, $ then $\, A
(i ,l)= A(i ,l') \, $ for $\, i \in \{ 1,\ldots ,u\} \, $ and $\,
l,l'\in \{ 1,\ldots ,\mu _{q_{i }} \} \, $; we denote this matrix
by $ \, A(i )\, $. \newline We will denote by $ \cSC _B$ the
subspace of all $A\in \cC _{B } $ such that $A(i )$ is lower
triangular for $i =1,\ldots ,u$. Moreover we will set
$$ \cSN _B=\cSC _B\; \cap \; \cN _B\
.$$ \begin{example}\label{ex4} {\rm If $B=(3,3,3,2)$ we have that $A\in \cN _B$
iff there exists a set  $\{ a_{h,k}^l\in K\  |\ (h,k,l)\in \{
1,2,3,4\} ^2 \times \{ 1,2,3\}\} $ such that $A$  is the
matrix:\vspace{2mm}
$${\small {\pmatrix{ a_{11}^1 & a_{11}^2 & a_{11}^3 & | &  a_{12}^1 & a_{12}^2
& a_{12}^3 & | &  a_{13}^1 & a_{13}^2 & a_{13}^3 & | & a_{14}^1 &
a_{14}^2\cr  & a_{11}^1 & a_{11}^2 & | &   & a_{12}^1 & a_{12}^2 &
| &
  & a_{13}^1 & a_{13}^2 & | &  & a_{14}^1 \cr & & a_{11}^1 & | & & & a_{12}^1 & | & &  & a_{13}^1 & | & &
\cr - & -& -& -& -& -& -& -& -& -& - & -& -& - \cr a_{21}^1 &
a_{21}^2 & a_{21}^3 & | &  a_{22}^1 & a_{22}^2 & a_{22}^3 & | &
a_{23}^1 & a_{23}^2 & a_{23}^3 & | & a_{24}^1 & a_{24}^2\cr &
a_{21}^1 & a_{21}^2 & | & & a_{22}^1 & a_{22}^2 & | &
 &  a_{23}^1 & a_{23}^2 & | &  & a_{24}^1 \cr & & a_{21}^1 & | & & & a_{22}^1 & | & & & a_{23}^1 & | &  &
 \cr - & -& -& -& -& -& -& -& -&
-& - & -& -& - \cr a_{31}^1 & a_{31}^2 & a_{31}^3 & | & a_{32}^1 &
a_{32}^2 & a_{32}^3 & | & a_{33}^1 & a_{33}^2 & a_{33}^3 & | &
a_{34}^1 & a_{34}^2\cr  & a_{31}^1 & a_{31}^2 & | &  & a_{32}^1 &
a_{32}^2 & | &
 &  a_{33}^1 & a_{33}^2 & | &  & a_{34}^1 \cr & & a_{31}^1 & | & & & a_{32}^1 & | & & & a_{33}^1 & | & &
\cr - & -& -& -& -& -& -& -& -& -& - & -& -& - \cr & a_{41}^1 &
a_{41}^2 & | &  & a_{42}^1 & a_{42}^2 & | & & a_{43}^1 & a_{43}^2
& | & 0 & a_{44}^2\cr & &   a_{41}^1 & | & & & a_{42}^1 & | &  & &
a_{43}^1 & | & & 0 \cr
 }}}$$ and $${\small \pmatrix{a_{11}^1 & a_{12}^1 & a_{13}^1 \cr a_{21}^1 & a_{22}^1
& a_{23}^1 \cr a_{31}^1 & a_{32}^1 & a_{33}^1\cr} }\in N(3,K)\ .$$
For each $A\in \cN_B$ there exists $\, G\in \cC _B\, $ and a set
$\{ \bar a_{h,k}^l\in K\  |\ (h,k,l)\in \{ 1,2,3,4\} ^2 \times \{
1,2,3\} \} $
 such
that det $G\neq 0$ and $ G^{-1}AG$ is the following element of
$\cSN _B$:\vspace{2mm}
 $${\small {\pmatrix{ 0 & \bar a_{11}^2 & \bar a_{11}^3 & | & 0 &
\bar a_{12}^2 & \bar a_{12}^3 & | & 0 & \bar a_{13}^2 & \bar
a_{13}^3& | & \bar a_{14}^1 & \bar a_{14}^2\cr  & 0 & \bar
a_{11}^2 & | & & 0 & \bar a_{12}^2 & | &
 &  0 & \bar a_{13}^2 & | &  & \bar a_{14}^1 \cr & & 0 & | & &  & 0 & | & &  & 0 & |  & & \cr
- & -& -& -& -& -& -& -& -& -& - & -& -& - \cr \bar a_{21}^1 &
\bar a_{21}^2 & \bar a_{21}^3 &  | & 0 & \bar a_{22}^2 & \bar
a_{22}^3 & | & 0 & \bar a_{23}^2 & \bar a_{23}^3 & | & \bar
a_{24}^1 & \bar a_{24}^2\cr  & \bar a_{21}^1 & \bar a_{21}^2 & | &
& 0 & \bar a_{22}^2 & | &
 &  0 & \bar a_{23}^2 & | &  & \bar a_{24}^1 \cr & & \bar a_{21}^1 & | & & & 0  & | & & & 0 & | & &
 \cr - & -& -& -& -& -& -& -& -&
-& - & -& -& - \cr \bar a_{31}^1 & \bar a_{31}^2 & \bar a_{31}^3 &
| & \bar a_{32}^1 & \bar a_{32}^2 & \bar a_{32}^3 & | & 0 & \bar
a_{33}^2 & \bar a_{33}^3& | & \bar a_{34}^1 & \bar a_{34}^2\cr &
\bar a_{31}^1 & \bar a_{31}^2 & | & & \bar a_{32}^1 & \bar
a_{32}^2 & | &
 &  0 & \bar a_{33}^2 & | &  & \bar a_{34}^1 \cr & & \bar a_{31}^1  & | & &  & \bar a_{32}^1  & | & &  & 0  & |  &  &
 \cr
 - & -& -& -& -& -& -& -& -&
-& - & -& -& - \cr
 & \bar a_{41}^1 & \bar a_{41}^2 & | &  & \bar a_{42}^1 & \bar a_{42}^2 & | & & \bar a_{43}^1
& \bar a_{43}^2 & | & 0 & a_{44}^2 \cr & &   \bar a_{41}^1 & | & &
& \bar a_{42}^1 & | & & & \bar a_{43}^1 & | &  & 0  \cr
 }}}$$}\end{example}
By lemma \ref{1} we get the following result.
\begin{corollary}\label{fromlemma2} The subvariety $\cSE
 _B$ ($\cSN _B$) has not empty intersection with the orbit of any
element of $\cE _B$ ($\cN _B$).\end{corollary}
\subsection{Upper
triangular form of $ \cSE _B$, $ \cSN _B$}\label{upper}  There is a bijection between $\Delta _B$ and the set
$$\{ (i,j,l)\ |\ i\in \{ 1,\ldots , u \}  , j\in \{ q_i-q_{i-1}, \ldots ,1\} , l\in \{ 1, \ldots , \mu _{q_i}\} \} \ ; $$ this set will be identified with $\Delta _B$ and with the set of the indices of the rows and the columns. It is ordered according to the following rule:
$(i,j,l)<(i',j',l')$ iff one of the following conditions
holds: \begin{itemize} \item[c$_1$)] $\ i<i'\ $ (that is $\ \mu
_{q_i}>\mu _{q_{i'}} $);\item[c$_2$)] $\ i=i'\ $ and $\ j>j'$;
\item[c$_3$)] $ \ i=i'\ $, $\ j=j'\ $ and $\ l>l'$. \end{itemize}
We observe that for $i\in \{ 1,\ldots ,u\} $ an element of $B$ is
equal to $\mu _{q_i}$ iff there exists $j\in \{ q_i-q_{i-1},\ldots
,1 \} $ such that it is the $ (q_i-j+1)-$th element of $B$.  If we
represent an element $X$ of $M(n,K)$ in the block form $
X=(X_{h,k})$, $\ h,k=1,\ldots ,t $ then the entry of $X\,
v_{\mu _{q_{i'}},j'}^{l'}$ with respect to $v_{\mu _{q_i},j}^l$, that is the entry of $X$ of indices $((i,j,l), (i',j',l'))$, is
the entry of the matrix $X_{q_i-j+1, q_{i'}-j'+1}$ which, in this
matrix, has indices $\big( \mu _{q_i}-l+1, \mu _{q_{i'}}-l'+1\big)
$. The next Lemma explains which square blocks of the matrix $X\in
\cSE _B$ are nilpotent (see the matrix $A$ of Example \ref{ex4}).
\begin{lemma}\label{1'} If $i,i'\in \{ 1,\ldots ,u\} \, $,
$ j\in \{ q_i-q_{i-1},\ldots ,1\} \,$ and $ j'\in \{
q_{i'}-q_{i'-1}, \ldots ,1\}  $
 the maximum rank of  $X_{q_i-j+1,q_{i'}-j'+1}$ for $\, X\in \cSE
_B$ ($ X\in \cSN _B $) is:
  \begin{itemize} \item[a)] $\ \mu
_{q_{i'}}\ \ $ if $\ \ i<i' \ $ ($\ \mu_ {q_i}
>\mu _{q_{i'}}\ $); \item[b)] $\ \mu _{q_i}-1\
\ $ if $\ \ i=i'\ \ $ and $\ \ j\leq j'\ $; \item[c)] $\ \mu
_{q_i}\ \ $ if $\ \ i>i'\ $ or $\; $ if $\ \ i=i'\ $ , $\ j>j'\
$.\end{itemize} \end{lemma} \pf The claim is a consequence of the
fact that if $X\in \cU _B$ then $X\in \cSE _B$ iff $X(i,l)$ is
strictly lower triangular for $i=1,\ldots ,u$ (if $X\in \cSN _B$
then $X(i,l)=X(i,l+1)$ for $l= 1,\ldots , \mu _{q_i}-1 $).
\hspace{4mm} $\square $ \vspace{2mm}
\newline If $X$ is an endomorphism of $K^n$ and $\Lambda $ is a
basis of $K^n$ we will denote by $R _{X, \Lambda }$ the relation
in the set of the elements of $\Lambda $ defined as follows: $\,
w'\ R _{X, \Lambda } \ w\, $ iff $\; X\, w'\, $ has nonzero entry
with respect to $w$. By the form of the matrices of $\cU _B$ and
lemma \ref{1'} we get the following description of the elements of
$\cSE _B$.
\begin{corollary}\label{cor'} There exists a not empty
open subset of $\, \cSE _B\, $ ($\, \cSN _B\, $) such that if $X$
belongs to it and $\ v_{\mu _{q_i},j}^{l}\ ,\ v_{\mu _{q_{i'}},
j'}^{l'}\in \Delta _B\ $ then $\, v_{\mu _{q_{i'}},j'}^{l'}\ R
_{X, \Delta _B} \ \ v_{\mu _{q_i},j}^l\ $, that is the entry of $X$ of indices $((i,j,l),(i',j',l'))$ is not $0$,  iff one of the following
conditions holds:
\begin{itemize} \item[$\iota _1 $)] $\, i<i'\, $  and $\,  \mu _{q_i}-l\leq \mu _{q_{i'}}-l'$; \item[$\iota  _2$)]
$ \, i=i'\, $, $\, j\geq j' \, $  and $\, l>l'$ (that is $\, \mu
_{q_i}-l<\mu _{q_{i'}}-l'\, $); \item[$\iota _3 $)] $\, i=i'\ $,
$\, j<j'\, $ and $\, l\geq l'$ (that is $\, \mu _{q_i}-l\leq \mu
_{q_{i'}}-l'\, $);
 \item[$\iota _4$)] $\ i>i'\ $
  and $\ l\geq l'\ $.\end{itemize} \end{corollary} The subalgebra $\cSE _B$ ($\,
\cSN _B$) is a maximal nilpotent subalgebra of $\cU _B$ ($\cC
_B$); hence there exists a basis of $K^n$ with respect to which
all the elements of $\cSE _B\, $ ($\cSN _B$)
 are  upper triangular. A basis with this property can be obtained just replacing the order of $\Delta _B$ with the following new
order: $ (i,j,l)\, \prec \, (i',j',l') $ if one of the following conditions holds:
\begin{itemize} \item[$e_1$)] $ \ \mu _{q_i}-l<\mu _{q_{i'}}-l'\ $; \item[$e_2$)] $\ \mu _{q_i}-l=\mu _{q_{i'}}-l'\
$ and $\ i<i'\ $ (hence $l>l'$); \item[$e_3$)] $\ \mu _{q_i}-l=\mu
_{q_{i'}}-l'\ $, $\ i=i'\ $ (hence $l=l'$) and $\ j>j'\ $.
\end{itemize} Let $\Delta _{B,\prec }$ be the basis of $K^n$ which has the same elements as $\Delta _B$ but in the order $\prec \, $.
By Corollary \ref{cor'}, the representation of all the
elements of $\cSE _B\, $ ($\, \cSN _B$) with respect to $\Delta
_{B,\prec }$ is upper triangular. \begin{example}\label{ex5} {\rm $\ $ Let $n=13$ and $B=(4,3^2,2,1)$. In this case
we have
$$\Delta _B=\{{v_{4,1}^4, v_{4,1}^3,
v_{4,1}^2, v_{4,1}^1}, {v_{3,2}^3, v_{3,2}^2,
v_{3,2}^1},{ v_{3,1}^3, v_{3,1}^2,
v_{3,1}^1},{ v_{2,1}^2, v_{2,1}^1},{
v_{1,1}^1}\} \ ;  $$ an element of $\cSN _B$, considered as an endomorphism,
 has the following matrix with respect to $\Delta _B$:\vspace{2mm} {\small{$$ \pmatrix{ 0 & a & b & c & | & p &
q & r & | & v & z & w & | & \alpha & \beta & | & \lambda \cr & 0 &
a & b & | & & p & q & | & & v & z & | &   & \alpha & | & 0 \cr & &
0 & a & | & & & p & | & & & v & | & &
 & | & \cr  & & & 0 & | & & & & | & & & & | & &  & | & &
\cr - & -& -& -& -& -& -& -& -& -& - & -& -& -& -& - \cr & s & t &
u & | & 0 & i & l & | & 0 & m & n & | & \rho & \sigma & | & \pi
\cr & & s & t & | & & 0 & i & | & & 0 &
 m & | &  & \rho & | & 0 \cr & & & s & | & & & 0 & | & & & 0
& | & & & \cr - & -& -& -& -& -& -& -& -& -& - & -& -& -& -& - \cr
& j & y & k & | & d & e & f & | & 0 & g & h & | & \xi & \zeta & |
& \theta \cr & & j & y & | & & d & e & | &  & 0 & g & | & & \xi &
| & 0 \cr & & & j & | & & & d & | & & & 0 & | & & & | & \cr - & -&
-& -& -& -& -& -& -& -& - & -& -& -& -& - \cr & { 0} & \delta &
\epsilon & | & & \eta & \nu & | & & \tau & \gamma & | & 0 & o & |
& \phi \cr & & { 0} & \delta & | & & & \eta & | & & &  \tau & | &
& 0 & | & \cr - & -& -& -& -& -& -& -& -& -& - & -& -& -& -& - \cr
& { 0} & { 0} & \omega & | & & { 0} & \psi & | & & { 0} & \iota &
| & & \chi & | & 0 \cr }\ .\vspace{2mm}$$}} If we instead consider the basis $$ \Delta _{B,\prec }=\{
{v_{4,1}^4, v_{3,1}^3, v_{3,2}^3, v_{2,1}^2,
v_{1,1}^1},{ v_{4,1}^3, v_{3,1}^2, v_{3,2}^2,
v_{2,1}^1},{ v_{4,1}^2, v_{3,1}^1, v_{3,2}^1},
{v_{4,1}^1}\}$$ then the previous endomorphism with
respect to $\Delta _{B,\prec } $ has the following matrix:\vspace{2mm} {\small {$$\pmatrix{0
& v & p & \alpha & \lambda & | & a & z & q & \beta & | & b & w & r
& | & c \cr & 0 & d & \xi & \theta & | & j & g & e & \zeta & | & y
& h & f & | & k \cr & & 0 & \rho & \pi &  | & s & m & i & \sigma &
| & t & n & l & | & u \cr & & & 0 & \phi & | & { 0} & \tau & \eta
& o & | & \delta  & \gamma  & \nu & | & \epsilon \cr  & & & & 0 &
| & { 0} & { 0} & { 0} & \chi & | & { 0} & \iota & \psi & | &
\omega \cr - & -& -& -& -& -& -& -& -& -& - & -& -& -& -& - \cr &
& & & & | & 0 & v & p & \alpha & | & a & z & q & | & b \cr & & & &
& | & & 0 & d & \xi & | & j & g & e & | & y \cr &  & & & & | & & &
0 & \rho & | & s & m & i & | &
 t \cr &  & & & & |&  & & & 0 & | & {0} & \tau & \eta &  | &
 \delta \cr - & -& -& -& -& -& -& -& -& -& - & -& -& -& -& - \cr &
 & & & &  | & & & & & | & 0 & v & p & | & a \cr &
& & & &  | & & & &  & | &  & 0 & d & | & j \cr &  & & & & | & & &
& & | & & & 0 & | & s \cr  - & -& -& -& -& -& -& -& -& -& - & -&
-& -& -& - \cr &  & & & & | & & & & & | & & & & | & 0 \cr  } \
.$$}}}\end{example}  In the following part of this Section and in the next one we will represent
any endomorphism of $K^n$ with respect to the basis $\Delta
_{B,\prec }\ $.  For $\ h=0,\ldots , \mu _{q_1}-1\ $ let
$$\Delta _{B,h}=\{ v_{\mu _{q_i},j}^l\in \Delta _B \ | \ \mu
_{q_i}-l=h\} $$ which will be considered as an ordered set with the order induced by $\prec $. We
set $t_h= \Big| \Delta
_{B, h}\Big| $ and $t_{-1}=0$; then $t_h\, \geq \, t_{h'}$ if $\,
h<h'$. Let $\pi _h$ be the canonical projection of $K^n$ onto
$\langle \Delta _{B, h}\rangle $; for $X\in M(n,K)\ $ and $\
h,k\in \{ 0, \ldots , \mu _{q_1}-1\} \ $ let $ \ X_{h,k}=\pi _h\,
\circ \, X|_{\langle \Delta _{B, k} \rangle }\ $; we consider $X$
as a block matrix:
$$X=(X_{h,k})\ , \ h,k\in \{ 0, \ldots , \mu _{q_1}-1\} \ .$$ In
the remainder of this Subsection we will describe the
representation of the algebra $\cSN _B$ with respect to $\Delta
_{B,\prec }\ $, as it appears in Example \ref{ex5}. This will also give a description of
 the representation of $\cSE _B$ with respect to $\Delta
_{B,\prec }\ $, since $\cSE _B$ is the minimal subspace of $\cN $ which contains $\cSN _B$ and is
defined by the condition that some coordinates are $0$.\newline
 Let $\, A\in \cSN _B\, $. If we
cancel the row and the column of index $(1,q_1,l) $
 for $ l=1,\ldots ,\mu _{q_1} $, that is the first
row of each row of blocks and the first column of each column of
blocks, we get a matrix $\, A^{\ast }\, $ of $\, \cSN _{B^{\ast }}\,
$, where $\, B^{\ast }=(\mu _2,\ldots ,\mu _t)\, $. If we
cancel the row and the column of index $(q_u,h+1,1)$
(that is of index ${\displaystyle \sum _{l=0}^{h}t_l}$ with respect to $\Delta _{B,\prec }$) for $h=0,\ldots ,\mu _{q_u}-1 $
we get a matrix $\,  A_{\ast}\, $ of $\, \cSN _{B_{\ast }}\, $, where
$\, B_{\ast }=(\mu _{1},\ldots ,\mu _{t-1})\, $.
\begin{proposition}\label{forma} With respect to $\Delta _{B,\prec }$ the
variety $\cSN _B$ is the affine space of all the strictly upper
triangular matrices $A$ such that:\begin{itemize} \item[i)] if $\,
\mu _{q_{i'}}-\mu _{q_i}>l'-l>0\, $ (hence $\mu _{q_{i'}}-\mu
_{q_i}>1$) then the entry of $\, A\, v_{\mu _{q_{i'}},j'}^{l'}\, $ with
respect to $\, v_{\mu _{q_i},j}^l\, $ is $0$; \item[ii)] for $\,
k\in \{ 1,\ldots ,\mu _{q_1}-1\} \, $ and $\, h\in \{ 1,\ldots
,k\} \, $ the entry of $\, A_{h,k}\, $ of indices $\, (i,j)\, $ is
equal to the entry of $\, A_{h-1,k-1}\, $ of indices $\, (i,j)\,
$.
\end{itemize} \end{proposition}
\pf The claim i) can be proved by Corollary \ref{cor'}; the claim
ii) can be proved by Lemma \ref{2}, Corollary \ref{cor'} and
induction on $n$ (we can consider $A^{\ast }$ or $A_{\ast }$).
\hspace{4mm} $\square $\vspace{2mm} \newline
The following result is a direct consequence of Proposition
\ref{forma}.
\begin{corollary}\label{description} If $\, X\in M(n,K)\, $ then $\, X\in
\cSN _B\, $ if and only if $\, X\, $ is a strictly upper triangular matrix
such that $\, X_{h+1,k+1}\, $ is a submatrix of $\, X_{h,k}\, $  and, if
$\, X_{h+1,k+1}\neq X_{h,k}\, $, then  $\, X_{h,k}\, $ is obtained by writing $0$'s
under $\, X_{h+1,k+1}\, $ and, if needed, any other columns on its right, for $\, h,k=0,\ldots
, \mu _{ q_1}-2\, $.\end{corollary}
We denote by $X=(x_{i,j})$, $i,j=1,\ldots ,n$ the generic element of $\cSE _B$ (with respect to $\Delta _{B,\prec }$),
which represents a morphism from an affine space to $\cN $ such that $x_{i,j}$ is $0$ or is just one of the coordinates of that affine space ($x_{i,j}\neq x_{i',j'}$ if $(i,j)\neq (i',j')$ and $x_{i,j}$ or $x_{i',j'} $ is not $0$). Instead we denote by $A=(a_{i,j})$, $i,j=1,\ldots ,n$  the generic
element of $\cSN _B$ (with respect to $\Delta _{B,\prec }$),
which represents a morphism from an affine space to $\cN $ such that $a_{i,j}$ is $0$ or is just one of the coordinates of that affine space. The following result is also a consequence of Proposition \ref{forma}.
\begin{corollary}\label{distance} Let $\, i,i',j,j'\in \{ 1,\ldots ,n\} \, $, $\, j\leq j'\, $. If $\, a_{i,j}=a_{i',j'} \, $ and $\, h\in \{ 0, \ldots ,\mu _{q_1}\} \, $ is such that  $\, j \in \Bigg\{ {\displaystyle \sum _{l=0}^{h-1}t_l}, \ldots , {\displaystyle \sum _{l=0}^{h}t_l} \Bigg\} \, $ then $\, j'-j=t_h\, $; replacing $a_{i,j}$ with $a_{j,i}$ and $a_{i',j'}$ with $a_{j',i'}$ in the previous assertion one gets an assertion which is still true. \end{corollary}
For $h=0,\ldots ,\mu _{q_1}-1$ let $B^{(h)}=(\mu
_1-h,\ldots ,\mu _t-h)$
 (omitting the values which are not positive). By Proposition
\ref{forma} we get the following result. \begin{lemma} \label{ind}
If $h\in \{ \mu _{q_1}-1,\ldots ,0\} $ the submatrix of $X$ (of $A$)
obtained by choosing the last $n- {\displaystyle \sum _{l=0}^{h-1} t_{l}}\, $ rows and columns is the
generic element of $\cSE _{B^{(h)}}$ (of $\cSN _{B^{(h)}}$).
\end{lemma}
By Corollary \ref{description} we get the following result.
\begin{proposition}\label{zeroblocks}
If $h\in \{ 1,\ldots , \mu _{q_1}-2\} $ and $ \, k\in \{h+1, \ldots ,\mu _{q_1}-1\} $ then $$A_{h,k}={\displaystyle \pmatrix{ A_{h,k}^{(\mu _{q_1})}  \cr \cr A_{h,k}^{(\mu _{q_1}-1)}
\cr \cr A_{h,k}^{(\mu _{q_1}-2)}\cr \cr \vdots \cr \cr  A_{h,k}^{(k+1)}\cr }}$$ (omitting the blocks different from $A_{h,k}^{(\mu _{q_1})}$ if $k=\mu _{q_1}-1$),
where $A_{h,k}^{(\mu _{q_1})}$ has $t_{\mu _{q_1}-1-(k-h)}$ rows and no zero entries while for $\lambda = \mu _{q_1}-1, \ldots ,k+1$ the block $A_{h,k}^{(\lambda )}$ has $t_{\lambda -1-(k-h)}-t_{\lambda -(k-h)} $ rows, has the first $t_{\lambda }$ columns equal to $0$ and has the other entries different from $0$ (it is possible that $t_{\lambda -1-(k-h)}-t_{\lambda -(k-h)} $ is $0$). \end{proposition}
It is obvious that if $h\in \{ 1,\ldots , \mu _{q_1}-3\} $, $ \, k\in \{h+1, \ldots ,\mu _{q_1}-2\} $ and $\lambda \in \{ k+2, \ldots , \mu _{q_1}\} $ then $A_{h,k}^{(\lambda )}$ and $A_{h+1,k+1}^{(\lambda )}$ have the same number of rows and the same number of zero columns, while if $h\in \{ 1,\ldots , \mu _{q_1}-3\} $, $ \, k\in \{h+1, \ldots ,\mu _{q_1}-2\} $ and $\lambda \in \{ k+1, \ldots ,\mu _{q_1}-2\} $ then $A_{h,k}^{(\lambda )}$ and $A_{h,k+1}^{(\lambda +1)}$ have the same number of rows (but respectively $t_{\lambda }$ and $t_{\lambda +1}$ zero columns).\newline
For $h\in \{ 1,\ldots , \mu _{q_1}-2\} \, $, $ \, k\in \{h+1, \ldots ,\mu _{q_1}-1\} $  we will denote by $l (h,k) $ the maximum of the set of all $l\in \{ h, \ldots ,k-1\} $ such that $t_{l}\neq t_k$. By this definition and Proposition \ref{zeroblocks} we get the following result.
\begin{corollary} \label{projection} If $h\in \{ 1,\ldots , \mu _{q_1}-2\} \, $, $ \, k\in \{h+2, \ldots ,\mu _{q_1}-1\} \, $ and $\, \lambda \in \{  k+1, \ldots ,\mu _{q_1}-1\, \} $ then all the blocks $A_{l(h,k ) ,k'}^{(\lambda )}\, $, $k'=l (h,k )+1, \ldots , k-1$ are zero blocks.\end{corollary}
\begin{corollary}\label{cor''} If  $h,k\in \{ 0,\ldots ,\mu _{q_1}-1\}
$ and the entry of $A_{h,k}$ of indices $(i,j)$
is $0$ for all $A\in \cSN _B$ then, for all $A\in \cSN
_B$: \begin{itemize} \item[1)] the entry of $
A_{h,k} $ of indices $ (i',j') $ is $0$ for $ i'=i,\ldots
,t_h $ and $ j'=1,\ldots ,j$; \item[2)] the
entry of indices $ (i,j) $ of $A_{h,k-1} $ (if $ k\neq
0$) and of $A_{h+1,k} $ (if $h\neq \mu _{q_1}-1
$ and $ i\leq t _{h+1}$) is also $0$. \end{itemize}
\end{corollary} \pf We can prove the claim by
Proposition \ref{forma}, using for simplicity induction on $n$ and
lemma \ref{ind}. By ii) of Proposition \ref{forma} in order to
prove 2) it is enough to prove the following claim: if $ k\in
\{ 1,\ldots ,\mu _{q_1}-1\} $ and the entry of indices $
(i,j)$ of $A_{0,k}$ is $ 0 $ for all $A\in \cSN _B
$ then the entry of indices $(i,j)$ of $A_{0,k-1}$ is
also $0 $ for all $A\in \cSN _B$. It can be proved by
induction on $t _{ 0}-i $ (if $ i=t _{0} $
the claim is true, since the claim is true for $A_{(1)}$ by the inductive hypothesis). \hspace{4mm} $\square $
\begin{corollary}\label{corresp} If $h \in \{ 1,\ldots ,\mu _{q_1}-1\} $, $\, k\in \{ 0,\ldots ,h\} \, $ and $$i\in \Bigg\{ {\displaystyle \sum_{l=0}^{k-1}t_l} , \ldots , {\displaystyle \sum_{l=0}^{k-1}t_l}+t_h\Bigg\} \ , \quad j\in \Bigg\{ {\displaystyle \sum_{l=0}^{h-1}t_l}, \ldots , {\displaystyle \sum_{l=0}^{h}t_l}\Bigg\} $$ then the entry of $A$ of indices $(i,j)$ is not $0$. \end{corollary}
\pf By Corollary \ref{description} the claim is obviously true for small values of $n$ (of $\mu _{q_1}$), hence we prove it by induction on $n$. By the inductive hypothesis the claim is true for the generic element $A_{\ast }$ of $\cSN _{B_{\ast }}$. Then by Corollary \ref{description} we get the claim also for $A$. \hspace{4mm} $\square $
\begin{corollary}
\label{added} For $i=1,\ldots ,n $ the entry of $A$ of indices
$(i,i+1)$ is not $0$ iff $\, i\neq  {\displaystyle \sum _{l=0}^ht_l }\, $ for $h=0,
\ldots , \mu _{q_1}-1$ or $\, {\displaystyle i=\sum _{l=0}^ht_l }\, $ and $t_h\leq t_{\mu _{q_1}-2}$ (hence the submatrix of $A$ obtained by choosing the rows and columns of indices greater than ${\displaystyle \sum _{l=0}^{\breve h}} t_l$ has kernel of dimension $1$).\end{corollary}
\subsection{Not empty intersection between $ \cSN _B$ and the maximum nilpotent orbit which intersects $ \cSE _B$}\label{indepToep}
In the set $\{ 1,\ldots ,n\} ^2$ we will consider the following strict partial order: $(i,j)<(i',j')$ iff $i\geq i'$,  $\, j\leq j'$ and $(i,j)\neq (i',j')$. It induces a relation in the set of the coordinates of $\cSN _B$, defined as follows: if $a$ and $a'$ are coordinates of $\cSN _B$, we will say that $a$ is smaller than $a'$ if there exist $(i,j), \ (i',j')\, \in \ \{ 1,\ldots ,n\} ^2$ such that $(i,j)<(i',j')$, $a$ is the entry of $A$ of indices $(i,j)$ and $a'$ is the entry of $A$ of indices $(i',j')$. We observe that if $a$ is smaller than $a'$ and $a'$ is both the entry of indices $(i',j')$ and the entry of indices $(\iota ', \gamma ' )$ with $\iota '< i'$ then there exists $(\iota , \gamma )< (\iota ', \gamma ' )$ such that the entry of $A$ of indices $(\iota , \gamma )$ is $a$; this is a consequence of Proposition \ref{forma}. It is obvious that if $a$ is smaller than $a'$ then $a'$ is not smaller than $a$. Let $a,\ a'$ be the entries of $A$ of indices $(i,j)$ and $(i',j')$ with $(i,j)<(i',j')$; let $a',\  a''$ be the entries of $A$ of indices $(\iota ', \gamma ')$ and $(i'', j'')$ with $(\iota ', \gamma ')<(i'', j'')$ and $\iota '\neq i'$. If $\iota '<i'$ then by the previous observation we get that $a$ is smaller than $a''$; if $\iota '>i'$ and $a$ is not smaller than $a''$ then $(i'',j'')$ is an entry of the first row of blocks of $A$ and we would get that $a$ is smaller than $a''$ by adding another row and column of blocks as a first row and column, according to Corollary \ref{description}. Hence this relation in the set of the coordinates of $\cSN _B$ can be extended to a partial order. If $a$ is not smaller than $a'$ and $a\neq a'$ we will say that $a$ is greater than $a'$.\newline
The following result is a consequence of Corollary \ref{cor''}.
\begin{corollary} \label{trasc} If $I\subseteq \{ 1,\ldots ,n-1\} $,
$J\subseteq \{ 2,\ldots ,n\} $  and $\overline i$, $\overline j$ are respectively the minimum element of $I$ and the maximum element of $J$ then: \begin{itemize}  \item[a)] if there exists $(\overline i,k)\in I\times J $ such
that $a_{\overline i,k}\neq 0$ then there exists $j\in J $ such that $a_{\overline i,
j}$ is different from all the
entries of the set $$\{ a_{i',j'}\ | \  (i',j')\in I\times J,\
(i',j')\neq (\overline i, j)\} \ ; $$
\item[b)] if there exists $(k,\overline j)\in I\times J$ such
that $a_{k,\overline j}\neq 0$ then there exists $i\in I$ such that $a_{i,
\overline j}$ is different from all the
entries of the set $$\{ a_{i',j'}\ | \  (i',j')\in I\times J,\
(i',j')\neq (i, \overline j)\} \ . $$\end{itemize}
\end{corollary}
\pf We first prove a). Let us consider the element $h_k\in \{ \mu _{q_1} -1, \ldots ,0\} $ such that
${\displaystyle \sum _{l=0}^{h_k-1}t_l<k\leq \sum _{l=0}^{h_k}t_l}$;  let $h_J$ be the maximum of the set $$
\Big\{ 0\Big\} \, \cup \, \Bigg \{ h\in \{ \mu _{q_1} -1, \ldots
,1\} \ | \ J\not \subseteq \Bigg\{ 1, \ldots , {\displaystyle \sum _{l=0}^{h-1}t_l}\, \Bigg\} \Bigg\} \
.$$ If $h_J=h_k\, $ by ii) of Proposition \ref{forma} we can set
$j=k$, since $a_{\overline i,k}$ is different from all the other entries of
the submatrix $(a_{i',j'})$, $i'\in I$, $j'\in J$. Hence we can
prove the claim by induction on $h_J-h_k$. If $a_{\overline i,k}$ is not
different from all the entries of the set $$\{ a_{i',j'}\ | \  i'\in I
, \ j'\in J,\ (i',j')\neq (\overline i, k)\} $$ then by ii) of
Proposition \ref{forma} there exist $k'\in J\, $ and $\, h_{k'}\in
\{ \mu _{q_1}-1,\ldots ,0\} \, $ such that $\, h_{k'}>h_k\, $, $\;
{\displaystyle \sum _{l=0}^{h_{k'}-1}t_l}<k'\leq {\displaystyle \sum _{l=0}^{h_{k'}}t_l}\, $, $\;
k'-{\displaystyle \sum _{l=0}^{h_{k'}-1}t_l}=k-{\displaystyle \sum _{l=0}^{h_{k}-1}t_l}$. But then by 2) of Corollary
\ref{cor''} we get that $a_{\overline i,k'}\neq 0$. Hence the claim follows
by the inductive hypothesis. The claim b) can be proved in the same way.\
 \hspace{4mm} $\square $\vspace{2mm}\newline
We define the following equivalence relation in the set $\{ 1,\ldots ,n\} $: we say that $i$ and $j$ are similar if there exist an entry of $A$ of column index $i$ and an entry of $A$ of column index $j$ which are equal; by Corollary \ref{description} one gets an equivalent condition by replacing the word "column" with the word "row". By Corollary \ref{trasc} we get the following result.
\begin{corollary}\label{trasc'} If $I\subseteq \{ 1,\ldots ,n-1\} $,
$\; J\subseteq \{ 2,\ldots ,n\} $ and $\; I'=\{ i_1,\ldots ,i_m\} $, $\; J'=\{ j_1,\ldots ,j_{m'}\} $ are respectively a maximal subset of $I$ of similar rows and a maximal subset of $J$ of similar columns, then: \begin{itemize} \item[a)] if $i$ and $j$ are respectively the minimum possible element of $I'$ and the maximum possible element of $J$ such that $a_{i,j}\neq 0$ then $a_{i,
 j}$ is different from all the
entries of the set $$\{ a_{i',j'}\ | \  (i',j')\in I\times J,\
(i',j')\neq (i,  j)\} \ ;  $$ \item[b)] the same claim as in a) holds if $i$ and $j$ are respectively the minimum possible element of $I$ and the maximum possible element of $J'$. \end{itemize}\end{corollary}
By Corollary \ref{trasc'} we get the following result.
\begin{corollary}\label{T}
The rational function $A$ over $M(n,K)$ has the property $\ast )$ defined in Subsection \ref{P}.
\end{corollary}
{\bf  Proof of Theorem \ref{indToep} } \
By Corollary \ref{T}, Proposition \ref{rank} and Theorem \ref{equalorbit}
 the maximum nilpotent orbit which intersects $\cSN _B$ is the same as the maximum nilpotent orbit which intersects $\cSE _B$, which is the statement of Theorem \ref{indToep} . \hspace{4mm} $\square $
\section{Proof of Theorem \ref{eend} }\label{conjecture}
\subsection{The graph associated to $B$}
Let $R_B$ be the relation in the set of the elements of $\Delta
_B$ defined as follows:
$$v_{\mu _{q_{i'}}, j'}^{l'}\ R _B \ v_{\mu
_{q_i},j}^l \ \Longleftrightarrow \   \iota _1) \ \mbox{\rm or }
\iota _2) \ \mbox{\rm or } \iota _3) \ \mbox{\rm or } \iota _4)\
\mbox{\rm of Corollary \ref{cor'} holds} \ .
$$ \begin{proposition} \label{order} The relation $R _B$ in the
set of the elements of $\Delta _B$ is a strict partial order.
\end{proposition} \pf The relation $R _B $ is obviously
antisymmetric; the condition $\iota _1)$ implies $l>l'$, hence it
is also transitive. \hspace{4mm} $\square $\vspace{2mm} \newline
The relation $R_B$ describes a generic element of $\cSE_B$, hence
by Theorem \ref{indToep} we get the following result.
\begin{corollary}\label{rel} The maximum nilpotent orbit of the
elements of $\cSN _B$ is determined by the relation $R _B$.
\end{corollary} We will write the vertices of the graph of $R_B$
forming a table in the following way: the
elements of $\N \cup \{ 0\} $ will be the indices of the rows and $\mu _{q_u},\ldots ,\mu _{q_1}$ will be
the indices of the columns; the element $v_{\mu _{q_i},j}^l$ of $\Delta _B$ will be
written in the $\mu _{q_i}-$th column and in the row whose index is
the maximum number $m$ such that there exist elements of $\Delta
_B$ whose images under $X^m$ have nonzero entry with respect to
$v_{\mu _{q_i},j}^l$ for some $X\in \cSE _B$. The
 graph of $R _B$ could be obtained by writing arrows on this table
 according to Corollary \ref{cor'}. We will say that $v_{\mu _{q_{i'}}, j'}^{l'}$
 precedes $v_{\mu
_{q_{i}}, j}^{l}$ in the graph of $R_B$ if the index of the row of $v_{\mu
_{q_{i'}}, j'}^{l'}$ in the previous table is less than the index of the row of $v_{\mu _{q_{i}},
j}^{l}$.  If there is an arrow from $v_{\mu _{q_{i'}},
j'}^{l'}$ to $v_{\mu _{q_{i}}, j}^{l}$ then $v_{\mu _{q_{i'}}, j'}^{l'}$
 precedes $v_{\mu
_{q_{i}}, j}^{l}$. \begin{lemma}
\label{graph} The vectors of $\Delta _B$ appear in the graph
of $R_B$ according to the following rules: \begin{itemize}
\item[a)] $v_{\mu _{q_{i}}, j'}^{l'}$ precedes $v_{\mu _{q_{i}},
j}^{l}$ if $l'<l$ or if $l'=l$ and $j'>j$;
\item[b)] if $\mu _{q_{i'}}<\mu _{q_i}$ and $\mu _{q_{i'}}-l'\geq
\mu _{q_i}-l$ then $v_{\mu _{q_{i'}}, j'}^{l'}$ precedes $v_{\mu
_{q_{i}}, j}^{l}$; \item[c)] if $\mu
_{q_{i'}}>\mu _{q_i}$ and $l'\leq l$ then $v_{\mu _{q_{i'}},
j'}^{l'}$ precedes $v_{\mu _{q_{i}}, j}^{l}$.
\end{itemize} \end{lemma}\pf The claim a) follows by $\iota _2$)
and $\iota _3$) of Corollary \ref{cor'}; the claims b) and c)
follow respectively from $\iota _1$) and $\iota _3$) of the same
Corollary.\hspace{5mm} $\square $ \vspace{2mm}\newline {\bf
Example 6}\hspace{4mm} For  $\, B=(7,5,2)\, $, $\, B=(2^2,1)\, $,  $\,
B=(4,2^2,1)\, $, $\, B=(2^3)\,
$, $\, B=(5,2^3) \,$,  $\, B=(6,5,2^3)\, $, $\, B=(3^2,2,1)\, $,
$\, B=(8^2,6^4, 3^2,2,1)\, $,  $\, B=(4,3^2,2,1)\, $,
 $\,
B=(5,4,3^2,2,1)\, $ and $\, B=(17,15,13,5,4,3^2,2,1)\, $ we respectively get the following
graphs (where we omit the arrows and we put $\circ   $ in front of the elements of $\Delta ^{\circ }_B$):
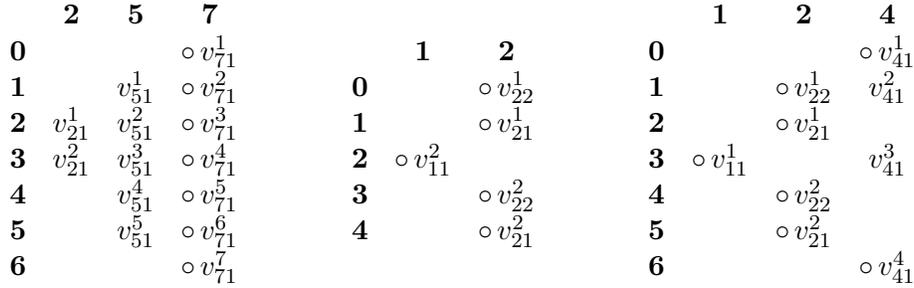
$$ \begin{array}{cccc} & \mbox{\bf 2} & \mbox{\bf 5}& \mbox{\bf 7}
\\ \mbox{\bf 0} & & &  \circ \, v_{71}^1\\ \mbox{\bf 1} & & v_{51}^1 &
\circ \, v_{71}^2  \\ \mbox{\bf 2} & v_{21}^1  & v_{51}^2 & \circ \, v_{71}^3\\
\mbox{\bf 3} & v_{21}^2 & v_{51}^3 & \circ \, v_{71}^4 \\ \mbox{\bf 4} &
& v_{51}^4 & \circ \, v_{71}^5\\ \mbox{\bf 5} &  & v_{51}^5 & \circ \, v_{71}^6\\
\mbox{\bf 6} & & & \circ \, v_{71}^7 \end{array}\qquad \quad
\begin{array}{ccc} & \mbox{\bf 1} & \mbox{\bf 2} \\ \mbox{\bf 0} &
&\circ \, v_{22}^1
\\  \mbox{\bf 1} &  & \circ \, v_{21}^1
\\ \mbox{\bf 2} & \circ \, v_{11}^2&  \\ \mbox{\bf 3} &  & \circ \, v_{22}^2 \\ \mbox{\bf 4} &
&\circ \, v_{21}^2
 \end{array} \qquad \quad  \begin{array}{cccc} & \mbox{\bf 1} & \mbox{\bf 2}
& \mbox{\bf 4}\\ \mbox{\bf 0} &  &  & \circ \, v_{41}^1\\ \mbox{\bf 1} &
 & \circ \, v_{22}^1 & v_{41}^2
\\ \mbox{\bf 2} & & \circ \, v_{21}^1 & \\ \mbox{\bf 3} & \circ \, v_{11}^1& &
v_{41}^3  \\ \mbox{\bf 4} & & \circ \, v_{22}^2 & \\ \mbox{\bf 5} & & \circ \, v_{21}^2 & \\
\mbox{\bf 6} &  & & \circ \, v_{41}^4
\end{array} \vspace{2mm} $$
$$ \begin{array}{c}\begin{array}{cc} & \mbox{\bf 2}\\ \mbox{\bf 0} & \circ \, v_{23}^1 \\ \mbox{\bf
1}& \circ \, v_{22}^1 \\ \mbox{\bf 2} & \circ \, v_{21}^1 \\ \mbox{\bf 3} & \circ \, v_{23}^2
 \\ \mbox{\bf 4} & \circ \, v_{22}^2 \\ \mbox{\bf 5} &
 \circ \, v_{21}^2\end{array}\qquad
\begin{array}{ccc} & \mbox{\bf 2}& \mbox{\bf 5}\\ \mbox{\bf 0} & &\circ \, v_{51}^1  \\ \mbox{\bf
1}&\circ \, v_{23}^1 & v_{51}^2 \\ \mbox{\bf 2} &\circ \, v_{22}^1  &v_{51}^3  \\
\mbox{\bf 3} &\circ \, v_{21}^1 &
 \\ \mbox{\bf 4} &\circ \,  v_{23}^2 &v_{51}^4  \\ \mbox{\bf 5} & \circ \, v_{22}^2 & \\ \mbox{\bf 6} &\circ \, v_{21}^2 &  \\
 \mbox{\bf 7} &  & \circ \, v_{51}^5\end{array} \\ \\ \\ \\ \begin{array}{cccc}  & \mbox{\bf 2} & \mbox{\bf 5} & \mbox{\bf 6} \\
\mbox{\bf 0} & & &\circ \, v_{61}^1   \\
\mbox{\bf 1} & & \circ \, v_{51}^1  &  \\
\mbox{\bf 2} & v_{23}^1 & & \circ \, v_{61}^2 \\
\mbox{\bf 3} & v_{22}^1 & \circ \, v_{51}^2  &  \\
\mbox{\bf 4} & v_{21}^1 & &  \circ \, v_{61}^3  \\
\mbox{\bf 5} & v_{23}^2  & \circ \, v_{51}^3 &   \\
\mbox{\bf 6} & v_{22}^2 &  & \circ \, v_{61}^4
\\ \mbox{\bf 7} &v_{21}^2  & \circ \, v_{51} ^4   &  \\
\mbox{\bf 8} &  &  & \circ \, v_{61}^5 \\ \mbox{\bf 9} &    &\circ \, v_{51}^5  & \\
\mbox{\bf 10} &  & &\circ \,  v_{61}^6
\end{array}\\ \\ \\ \\ \begin{array}{cccc} & \mbox{\bf 1} & \mbox{\bf 2} &
\mbox{\bf 3}\\ \mbox{\bf 0} &  &  & \circ \, v_{32}^1 \\ \mbox{\bf 1} &
 & & \circ \, v_{31}^1
\\ \mbox{\bf 2} & & \circ \, v_{21}^1 & \\ \mbox{\bf 3} & v_{11}^1  & & \circ \, v_{32}^2
\\ \mbox{\bf 4} & & &\circ \,  v_{31}^2 \\ \mbox{\bf 5} & &
\circ \, v_{21}^2 & \\ \mbox{\bf 6} &  & & \circ \, v_{32}^3 \\ \mbox{\bf 7} &
 & & \circ \, v_{31}^3
\end{array}\end{array} \qquad \quad  \begin{array}{cccccc}  & \mbox{\bf 1} & \mbox{\bf 2} & \mbox{\bf 3} & \mbox{\bf 6} & \mbox{\bf 8}\\
\mbox{\bf 0} &  & & & & \circ \, v_{82}^1\\
\mbox{\bf 1} &  & & & & \circ \, v_{81}^1 \\
\mbox{\bf 2} &  & & & \circ \, v_{64}^1 & v_{82}^2 \\
\mbox{\bf 3} &  & & & \circ \, v_{63}^1 & v_{81}^2\\
\mbox{\bf 4} &  & & & \circ \, v_{62}^1 & \\
\mbox{\bf 5} &  & & & \circ \, v_{61}^1 & \\
\mbox{\bf 6} &  & & v_{32}^1  & \circ \, v_{64}^2 & v_{82}^3\\
\mbox{\bf 7} &  &  & v_{31}^1 & \circ \, v_{63}^2 & v_{81}^3 \\
\mbox{\bf 8} &  & v_{21}^1 &  & \circ \, v_{62}^2 &
\\ \mbox{\bf 9} &v_{11}^1  & &  & \circ \,  v_{61} ^2 & \\
\mbox{\bf 10} &  & &v_{32}^2  & \circ \, v_{64}^3 & v_{82}^4\\ \mbox{\bf 11} & &  & v_{31}^2 & \circ \, v_{63}^3 & v_{81}^4\\
\mbox{\bf 12} &  & v_{21}^2 &  &\circ \, v_{62}^3 & \\ \mbox{\bf 13} &  &   & & \circ \, v_{61}^3& \\
\mbox{\bf 14} &
&  &v_{32}^3 &\circ \, v_{64} ^4 & v_{82}^5\\ \mbox{\bf 15} &  &  &v_{31}^3 &\circ \, v_{63}^4& v_{81}^5\\ \mbox{\bf 16} &  & & & \circ \, v_{62}^4 & \\
\mbox{\bf 17} &  & & & \circ \, v_{61}^4 & \\ \mbox{\bf 18} &  &  & &
\circ \, v_{64}^5 & v_{82}^6
\\ \mbox{\bf 19} &  & & &\circ \,  v_{63}^5  & v_{81}^6\\ \mbox{\bf 20} & & & & \circ \, v_{62}^5 &  \\
\mbox{\bf 21} &  & & &\circ \, v_{61}^5 &  \\ \mbox{\bf 22} & & & &
\circ \, v_{64}^6 & v_{82}^7
\\ \mbox{\bf 23} &  & & & \circ \, v_{63}^6 & v_{81}^7  \\ \mbox{\bf 24} &  & & & \circ \, v_{62}^6 & \\
\mbox{\bf 25} & & & & \circ \, v_{61}^6 & \\ \mbox{\bf 26} &  & & & &
\circ \, v_{82}^8\\ \mbox{\bf 27} &  & & & & \circ \, v_{81}^8
\end{array}$$
$$\begin{array}{ccccc} & \mbox{\bf 1}& \mbox{\bf 2} & \mbox{\bf 3} & \mbox{\bf 4}\\
\mbox{\bf 0}&  & & &\circ \,  v_{41}^1 \\ \mbox{\bf 1}& & & \circ \, v_{32}^1  & \\ \mbox{\bf 2}& & & \circ \, v_{31}^1 & \\
\mbox{\bf 3}& &v_{21}^1 & & \circ \, v_{41}^2\\\mbox{\bf 4} & v_{11}^1 & &  \circ \, v_{32}^2 & \\
\mbox{\bf 5}
& & & \circ \, v_{31}^2 & \\ \mbox{\bf 6} & &v_{21}^2& &\circ \, v_{41}^3 \\ \mbox{\bf 7} & & &\circ \, v_{32}^3 & \\
\mbox{\bf 8} & & & \circ \, v_{31}^3& \\ \mbox{\bf 9} & & & &\circ \, v_{41}^4
\end{array}\qquad \qquad \quad \begin{array}{cccccc}  & \mbox{\bf 1}& \mbox{\bf 2}&\mbox{\bf 3} & \mbox{\bf 4} & \mbox{\bf 5}\\
\mbox{\bf 0} &  & & & &\circ \, v_{51}^1 \\
\mbox{\bf 1} & & & &\circ \, v_{41}^1 & \\ \mbox{\bf 2} & & & \circ \, v_{32}^1 & &  v_{51} ^2\\ \mbox{\bf 3} & & &\circ \, v_{31}^1 & & \\
\mbox{\bf 4} & & v_{21}^1& & \circ \, v_{41}^2& \\ \mbox{\bf 5} & v_{11}^1 & & \circ \, v_{32}^2& &v_{51}^3\\
\mbox{\bf 6} &
& & \circ \, v_{31}^2& & \\ \mbox{\bf 7} & & v_{21}^2& &\circ \, v_{41}^3 & \\ \mbox{\bf 8} &  & &\circ \,  v_{32}^3& & v_{51}^4 \\
\mbox{\bf 9} & & & \circ \, v_{31}^3& & \\ \mbox{\bf 10} & &  & & \circ \, v_{41}^4&
\\ \mbox{\bf 11} &  & & & &\circ \, v_{51}^5
\vspace{2mm}\end{array}$$ $$\begin{array}{ccccccccc}  & \mbox{\bf 1}& \mbox{\bf 2}&\mbox{\bf 3} & \mbox{\bf 4} & \mbox{\bf 5}& \mbox{\bf 13} & \mbox{\bf 15} & \mbox{\bf 17}\\
\mbox{\bf 0} &  & & & & & & & \circ \, v_{171}^1 \\
\mbox{\bf 1} &  & & & & & & \circ \, v_{151}^1 & v_{171}^2 \\
\mbox{\bf 2} &  & & & & & \circ \, v_{131}^1 & v_{151}^2 & v_{171}^3 \\
\mbox{\bf 3} &  & & & & \circ \, v_{51}^1 & v_{131}^2 & v_{151}^3 & v_{171}^4\\
\mbox{\bf 4} & & & &\circ \, v_{41}^1 & & v_{131}^3 & v_{151}^4 &
v_{171}^5\\ \mbox{\bf 5} & & & \circ \, v_{32}^1 & &  v_{51} ^2
& v_{131}^4 & v_{151}^5 & v_{171}^6\\ \mbox{\bf 6} & & &\circ \, v_{31}^1 & & & v_{131}^5 & v_{151}^6 & v_{171}^7\\
\mbox{\bf 7} & & v_{21}^1& & \circ \, v_{41}^2& & v_{131}^6 & v_{151}^7 &
v_{171}^8\\ \mbox{\bf 8} & v_{11}^1 & & \circ \, v_{32}^2& &v_{51}^3
& v_{131}^{7} & v_{151}^8 & v_{171}^{9}\\
\mbox{\bf 9} &
& & \circ \, v_{31}^2& & & v_{131}^8 & v_{151}^{9} & v_{171}^{10}\\ \mbox{\bf 10} & & v_{21}^2& &\circ \, v_{41}^3 & & v_{131}^{9} & v_{151}^{10} & v_{171}^{11}\\
\mbox{\bf 11} &  & & \circ \, v_{32}^3& & v_{51}^4 & v_{131}^{10} & v_{151}^{11} & v_{171}^{12}\\
\mbox{\bf 12} & & & \circ \, v_{31}^3& & & v_{131}^{11} & v_{151}^{12} &
v_{171}^{13}\\ \mbox{\bf 13} & &  & & \circ \, v_{41}^4 & & v_{131}^{12} &
v_{151}^{13} & v_{171}^{14}
\\ \mbox{\bf 14} &  & & & & \circ \, v_{51}^5 & & v_{151}^{14} &
v_{171}^{15} \\ \mbox{\bf 15} &  & & & & &\circ \, v_{131}^{13} & &
v_{171}^{16} \\ \mbox{\bf 16} &  & & & & & &  \circ \, v_{151}^{15} & \\
\mbox{\bf 17} &  & & & & & & & \circ \, v_{171}^{17}
\vspace{2mm}\end{array}$$
We recall the definitions of $\, (\tilde
i,\tilde {\epsilon })$, $\, \Delta _B^{\circ, 1}$, $\, \Delta
_B^{\circ, 2}$, $\, \Delta _B^{\circ, 3}$, $\, \Delta _B^{\circ
}$, $\, \widehat B$, $\, Q(B)=(\omega_1, \ldots ,\omega_z)\, $
which were given in section 1; the index of the last row of the
graph of $R_B $ where there are written some vectors is
$\omega_1-1$ and, by Theorem \ref{indToep}, the maximum partition
which is associated to elements of $\cSE _B$ is $Q(B)$.\newline By
the definition of the rows of the graph of $R_B$ any row of index
less than $\omega _1$ has at least one vector; by $\iota _2$) and
$\iota _3$) of Corollary \ref{cor'} the first vector and the last
vector of the $\mu _{q_i}-$th  column are respectively $v_{\mu
_{q_i},q_i-q_{i-1}}^1$ and  $v_{\mu _{q_i},1}^{\mu _{q_i}}$ for
$i=1,\ldots ,u$. Moreover we can observe the following less
obvious properties. We consider the injective map $\varphi \ :\
\Delta _B \longrightarrow \{ 0,\ldots ,\omega_1-1\} $ which
associates to any vector the index of its row in the graph of
$R_B$ and the injective map $\phi \, : \, \{ 0,\ldots ,\omega
_1-1\} \longrightarrow \Delta _ B$ such that $\phi (i)$ is the
first vector of the $i-$th row of the graph of $R_B$.
\begin{lemma} \label{next} The basis $\Delta _B$ with the order
$<$ has the following properties:
\begin{itemize} \item[i)] the restriction of the map $\varphi $ to
$\Delta _B^{\circ }$ is a bijection; \item[ii)] the restriction of
$\, \phi $ to $\{ 0,\ldots ,q_{\tilde i +\tilde{\epsilon }} -1\} $
preserves the order and its image is $\Delta _B^{\circ ,1}$; the
restriction of $\, \phi $ to $\{ \omega _1-q_{\tilde i
+\tilde{\epsilon }}+1,\ldots ,\omega _1\} $ reverses the order and
its image is $\Delta _B^{\circ,3}$; \item[iii)] if
 $\, i\in \{ 1,\ldots ,\tilde i+\tilde {\epsilon }\} $, $j\in \{ q_i-q_{i-1},\ldots ,1 \} $  and  $v_{\mu _{q_{i'}},j'}^{l'}$ is
written in the same row of the graph of $R_B$ as $\, v_{\mu
_{q_i},j}^{\mu _{q_i}}$ then $ i'\in \{1,\ldots ,i-1\} $ and
$l'\in \{ \mu _{q_i}+1, \ldots , \mu _{q_{i'}}-1\} $.\end{itemize}
\end{lemma} \pf We can write the graph of $R_B$ starting with writing
the vectors of $\Delta _B^{\circ }$ according to lemma
\ref{graph}. In fact, by Corollary \ref{cor'} and i) of
Proposition \ref{forma} if $\mu _{q_{i'}}-l'<\mu _{q_i}-l$ there
is no arrow from $v_{\mu _{q_{i'}},j'}^{l'}$ to $v_{\mu
_{q_{i}},j}^{l}$ and the same holds if $\mu _{q_{i'}}-\mu
_{q_i}>l'-l>0$; hence by the maximality of $2q_{\tilde i-1}+\mu
_{q_{\tilde i}}(q_{\tilde i}-q_{\tilde i-1})+\tilde {\epsilon }\,
\mu _{q_{\tilde i+1}}(q_{\tilde i+1}-q_{\tilde i})$ if the vectors
of $\Delta _B$ are written in the graph of $R_B$ according to
lemma \ref{graph} then each vector of $\Delta _B-\Delta_B^{\circ
}$ is in the same row as one of the vectors of $\Delta _B^{\circ
}$, which proves i). If $\mu _{q_i}<\mu _{q_{i'}}$ then $v_{\mu
_{q_{i'}},j'}^{1}$ precedes $v_{\mu _{q_i},j}^{1}$ in the graph of
$R_B$ by c) of lemma \ref{graph}, while $v_{\mu _{q_i},j}^{\mu
_{q_i}}$ precedes $v_{\mu _{q_{i'}},j'}^{\mu _{q_{i'}}}$ by b) of
lemma \ref{graph}, hence we get ii). Moreover for $i=1,\ldots
,\tilde i -1$, $j=q_i-q_{i-1},\ldots ,1$ and $l=\mu _{ q_i}-1,
\ldots ,1$ the index of the row of $v_{\mu _{q_i},j}^{l}$ is
less than the index of the row of $v_{\mu
_{q_{i+1}},q_{i+1}-q_i}^{\mu _{q_{i+1}}}$ by $\iota _2$), $\iota
_3$) of Corollary \ref{cor'} (see the graph of $R_B$ in the case
$B=(17,15,13, 5,4,3^2,2,1)$), hence we get iii).  \hspace{4mm}
$\square $\vspace{2mm} \newline By i) of lemma \ref{next} we get
Theorem \ref{O}, which was proved by Polona Oblak in \cite{Obl}.
\subsection{On the partition associated to elements of $\cE _B$ ($\cN _B$)}
Let $X$ be the matrix of an endomorphism of $K^n$ with respect to
the basis $\Delta _B$. We recall the following classical lemma (see \cite{Her}).
\begin{lemma} \label{H} If $\, X\in N(n,K)$ has index of nilpotency
$\omega $ and $v\in K^n $ is such that $X^{\omega -1}v\neq 0$
there exists a subspace $V$ of $K^n$ such that $X(V)\subseteq V$
and
$$V \; \oplus \; \langle v, Xv, \ldots , X^{\omega -1} v\rangle \, =
K^n\ .$$ \end{lemma} For $(X,v)\in \cSE  _B\, \times \, K^n$  let
$W_{X,v,0}=\{ 0\} $ and, for $i\in \N $, let $W_{X,v,i}=\langle
X^hJ^kv\, |\, h,k\in \N \, \cup \, \{ 0\} \, ,\, k<i \rangle $. The
subspace $W_{X,v,i}$ is stable with respect to $X$;  let $X_{v,i}$ be
the endomorphism of $\; {\displaystyle K^n/ W_{X,v,i}}\; $ defined
by $X_i(W_{X,v,i}+w)=W_{X,v,i}+Xw$.\newline
By Lemma \ref{H} we get the following result.
\begin{proposition} \label{j2} If $(X,v)\in \cSE _B\times K^n$ and $X^{\omega _1-1}v\neq
0$ the partition of $X$ is
 $(\omega_1, \nu _2, \nu _3, \ldots ,\nu _p)$ iff
the partition of $X_{v,1}$ is $(\nu _2, \nu _3, \ldots ,\nu _p)$.
\end{proposition}
 If $(X,v)\in N(n,K)\times K^n$
and $ K^n=\langle X^hJ^kv,\ h,k=1,\ldots ,n\rangle $ one says that
$v$ is cyclic for $X$.
\begin{proposition}\label{j} If
$(X,v)\in \cSE _B\, \times \, K^n$ and $v$ is cyclic for $X$ then:
\begin{itemize} \item[i)] $X$ has partition $(\nu_1,\ldots ,\nu_{p })$
iff the set $$\Omega _{(X,J)}=\{ \ X^hJ^kv\ |\ ,\ h= \nu _{k+1}-1,\ldots, 0\, , \, k=0,\ldots ,p-1 \ \} $$ is a Jordan basis for
$X$;  \item[ii)] if $X$ has partition $(\nu _1,\ldots ,\nu _{p
})$ and $h',k'$ are elements of $\N \, \cup \, \{ 0\} $ such that $k'\geq p $ or $h'\geq
\nu_{k'+1}$ then
$$ X^{h'}J^{k'}v\in \langle \ X^hJ^kv\, |\,
h=\nu _{k+1}-1,\ldots ,0\ ,\ k=0,\ldots ,k'-1\ \rangle \ .$$
\end{itemize} \end{proposition}
\pf  Since $v$ is cyclic for $X$ we have that
$${\displaystyle K^n/W_{X,v,i} }\, =\, \langle \ W_{X,v,i}+X^hJ^kv\ |\ h,k\in
\N \cup \{ 0\} \, ,\ k\geq i\ \rangle \ .$$ Let $X$ have partition
$(\nu _1,\ldots ,\nu _{p })$. Since $\nu _1$ is the index of
nilpotency of $X$ we have $X^{\nu _1-1}v\neq 0$, hence
$W_{X,v,1}=\langle v,Xv,\ldots ,X^{\nu_1-1}v\rangle $. Then $\nu_2$
is the index of nilpotency of $X_{v,1}$ by lemma \ref{H}. Hence
$X^{\nu _2-1}Jv\not \in W_{X,v,1}$. Similarly for $i=2,\ldots ,p -1$
we have that $\ \nu_ {i+1}$ is the index of nilpotency of $X_{v,i}$
and then $X^{\nu _{i+1}-1}J^iv\not \in W_{X,v,i}$.\hspace{4mm}
$\square $\vspace{2mm}\newline Let $K
_B$ be the subset of $K^n$ of all the vectors which have nonzero
entry with respect to $v_{ \mu _{q_1},q_1}^1$. The following Lemma is a
generalization of a known result for the elements of $\cN _B$.\begin{lemma}\label{j1} For all $v\in K _B$ the subset of
 $\cSE _B$ of all $X$ such
that $v$ is cyclic for $X$ is not empty.
\end{lemma} \pf In \cite{Na} it has been proved that the
subset of $\cSN _B\times K^n$ of all the pairs $(X,v)$ such that
$v$ is cyclic for $X$ is not empty. The projection of this subset
on $K^n$ is an open subset which is not empty, hence its
intersection with $K _B$ is not empty. Since any element of $K _B$
can be the $\mu_{q_1}-$th element of a Jordan basis for $J$ we get
that $K _B$ is contained in that projection. \hspace{4mm} $\square
$\vspace{2mm}\newline
We will consider pairs $(X,v)\in
N(n,K)\times K^n$ with the following property:
$$\begin{array}{l}   \mbox{\rm if $m\in \N\cup  \{ 0\} $,
$\, m\leq \omega_1-1$ then $ X^m\, v$ has nonzero entry with respect to} \\
\mbox{\rm the vector of  $\Delta _B^{\circ}$ which is written in
the $m-$th
 row of the graph of $R_B$} \ ;  \end{array}
$$ let $\cSE^{\star } _B$ be the open subset of $\, \cSE
_B\times K^n\, $ of all the pairs with this property.
\begin{proposition} \label{f1}
If $v\in K^n$ the fiber of $v$ with respect to the canonical
projection of $\cSE^{\star } _B$ into $K^n$ is not empty iff $v\in
K _B$.
\end{proposition} \pf By Corollary \ref{cor'} there exists a
not empty open subset of $\cSE _B$ ($\cSE^{\star } _B$) such that
if $X$ belongs to it then the vector $X\, v_{ \mu _{q_1},q_1}^1$ has
nonzero entry with respect to all the elements of $\Delta _B- \{
v_{ \mu _{q_1},q_1}^1\} $. Hence if $v\in K^n$ has nonzero entry
with respect to $v_{ \mu _{q_1},q_1}^1$ and $m,i\in \N \, \cup \,
\{ 0\} $, $ i\geq m$ then the condition that $X^m\, v$ has zero entry
with respect to a vector of the $i-$th row of the graph of $R_B$
isn't an identity with respect to $X$. \hspace{4mm} $\square
$
\begin{lemma}\label{directsum} If $\, (X,v)\in \cSE^{\star } _B$ then
$W_{X,v,1}\ \cap \ \langle \Delta _B-\Delta _B^{\circ }\rangle =\{
0\} $. \end{lemma} \pf Let $w$ be a non-zero element of $ W_{X,v,1}$;  let $m$ be the minimum element of $ \{
0,\ldots ,\omega_1-1\} $ such that $w$, if represented with respect to
$$\{ v, Xv, \ldots , X^{\omega _1 -1} v\} \ ,$$ has nonzero entry with
respect to $X^mv$.  By the definition of $\cSE ^{\star }_B$ and Corollary \ref{cor'} $w$
 has nonzero entry
with respect to the element of $\Delta _B^{\circ }$ which is
written in the row of index $m$; hence $w\not \in \langle \Delta
_B-\Delta _B^{\circ }\rangle$. \hspace{4mm} $\square $
\subsection{The maximum partition in $\cE _B$ ($\cN _B$)}
If $(X,v)\in \cSE^{\star } _B $ we will denote by $J_{X,v}$ the endomorphism of
${\displaystyle K^n/W_{X,v,1}}$ defined by $
J_{X,v}(W_{X,v,1}+w)=W_{X,v,1}+J(w)$ for $w\in \Delta _B-\Delta
_B^{\circ }$. We observe that if $X\in \cSN _B$ then $J_{X,v}$
commutes with $X_{v,1}$. In general $J_{X,v}$ is not in Jordan canonical
form with respect to the basis $\{ W_{X,v,1}+v_{\mu _{q_i},j}^l\ |\
v_{\mu _{q_i},j}^l\in \Delta _B -\Delta _B^{\circ }\} $; in fact
it is not true that $v_{\mu _{q_i},j}^{ \mu _{q_i}} \in W_{X,v,1} $
for $i = 1,\ldots ,\tilde i-1 $ and $j= q_i-q_{i-1}, \ldots ,1 $
(see the example with $B=(17,15,13, 5,4,3^2,2,1)$, where $\mu
_{q_{\tilde i}} =4$ and $\mu _{q_{\tilde i +\tilde {\epsilon
}}}=3$). In the next Corollary we will show that if $(X,v)\in \cSE^{\star } _B$
the partition of $n-\omega _1$ which corresponds to $J_{X,v}$ is the
partition $\widehat B$ which has been defined in the introduction.
\newline For $\, i=1,\ldots ,\tilde i-1\, $ and $\, j=q_i-q_{i-1},
\ldots ,1\, $ let $\, H(i,j)\, $ be the set of all $\,
(i',j',l')\, $ such that $\, i'\in \{ 1,\ldots i-1\} \, $, $\,
j'\in \{ q_{i'}-q_{i'-1}, \ldots ,1\} \, $ and $l'\in \{ 2, \ldots
, \mu _{q_{i'}}-\mu _{q_i}+1\} $.
\begin{lemma}   \label{nextt} If $(X,v)\in \cSE^{\star }_B$ then for $i=1,\ldots ,\tilde i-1$, $j=q_i-q_{i-1},\ldots ,1$
there exist $a(i,j,i',j',l')\in K$ for $(i',j',l')\in H(i,j)$ such
that:
$$ v_{\mu _{q_i},j}^{\mu _{q_i}}+\sum _{(i',j',l')\in H(i,j)}a(i,j,i',j',l') \, v_{\mu _{q_{i'}},j'}^{l'+\mu _{q_i}-1} \in W_{X,v,1}\ , $$
$$ v_{\mu _{q_i},j}^{\mu _{q_i}-1}+\sum _{(i',j',l')\in H(i,j)}a(i,j,i',j',l') \, v_{\mu _{q_{i'}},j'}^{l'+\mu _{q_i}-2} \not \in W_{X,v,1}\ . $$
\end{lemma} \pf Let $(X,v)\in \cSE _B^{\star }$. By iii) of lemma \ref{graph} for $i=1,\ldots ,\tilde i-1$ and $j=q_i-q_{i-1},\ldots ,1$
there exist $a\in K-\{ 0\} $ and $a(i,j,i',j',l')\in K$ for
$(i',j',l')\in H(i,j)$ such that $$X^{\varphi (v_{\mu
_{q_i},j}^{\mu _{q_i}})}\, v = a\, v_{\mu _{q_i},j}^{\mu
_{q_i}}+\sum _{(i',j',l')\in H(i,j)}a\cdot a(i,j,i',j',l')\,
v_{\mu _{q_{i'}},j'}^{l'+\mu _{q_i}-1}\ ;$$ this proves the first claim. For $m=0,\ldots, \varphi (v_{\mu
_{q_i},j}^{\mu _{q_i}-1})$ the vector
$X^{m}\, v $ has nonzero entry with respect to the vector of $\Delta _B^{\circ }$ which is written in the row of index $m$ of the graph of $R_B$, while it has zero entry with respect to the vectors of $\Delta _B^{\circ }$ which are written in the previous rows; hence $$v_{\mu _{q_i},j}^{\mu _{q_i}-1}\ \not \in \ W_{X,v,1}\,+\,  \Big\langle \,v_{\mu _{q_{i'}},j'}^{l'+\mu _{q_i}-2}\ , (i',j',l')\in H(i,j)\, \Big\rangle $$ which proves the second claim.
 \hspace{5mm} $\square $
\vspace{2mm} \newline In the following we will assume that for
$(X,v)\in \cSE _B^{\star }$, $\, i=1,\ldots ,\tilde i-1$, $j=q_i-q_{i-1},\ldots ,1$ and
$(i',j',l')\in H(i,j)$ the elements $a(i,j,i',j',l')\in K$ have
the property expressed in lemma \ref{nextt}.
\begin{corollary} \label{next'} If $(X,v)\in \cSE^{\star } _B$ and we set $$\hat v _{\mu
_{q_i},j}^l=\left\{
\begin{array}{l}v_{\mu _{q_i},j}^{l} + {\displaystyle \sum _{(i',j',l')\in H(i,j)}a(i,j,i',j',l')\, v_{\mu
_{q_{i'}},j'}^{l'+l-1}}\vspace{2mm}
 \\  \quad  \mbox{\rm if } i\in \{ 1,\ldots ,\tilde i-1\} ,\ j\in \{
q_i-q_{i-1},\ldots ,1\} , \ l\in \{ 1,\ldots ,\mu _{q_i}-1\} \\ \\
 v_{\mu _{q_i},j}^{l}  \vspace{1mm}\\  \quad  \mbox{\rm if } i\in \{ \tilde i+\tilde {\epsilon }+1,
\ldots ,u\} , \  j\in \{ q_i-q_{i-1},\ldots ,1\} , \ l\in \{
1,\ldots ,\mu _{q_i}\} \end{array} \right. $$ then the representation of
$J_{X,v}$ with respect to the basis
$$\widehat {\Delta }_{B,X,v}\, = \ \{ W_{X,v,1}+\hat v _{\mu _{q_i},j}^l\ |\ v_{\mu _{q_i},j}^l\in
\Delta _B-\Delta _B^{\circ }\} $$ (with the order induced by
$\Delta _B$) has Jordan canonical form and its partition is $\widehat B$.
\end{corollary} \pf  In lemma \ref{nextt} for $i=1,\ldots ,\tilde
i-1$ and $j=q_i-q_{i-1}, \ldots ,1$ we have proved the following
claims:
$$ J^{\mu _{q_i}-1}\, \Big( v_{\mu _{q_i},j}^{1}+\sum
_{(i',j',l')\in H(i,j)}a(i,j,i',j',l') \, v_{\mu
_{q_{i'}},j'}^{l'}\, \Big)  \in W_{X,v,1}\ , $$
$$ J^{\mu _{q_i}-2} \Big( v_{\mu _{q_i},j}^{1}+\sum _{(i',j',l')\in H(i,j)}a(i,j,i',j',l')\, v_{\mu _{q_{i'}},j'}^{l'} \, \Big)
\not \in W_{X,v,1}\ .
$$   By the definition of $J_{X,v}$
they imply that $(J_{X,v})^{\mu _{q_i}-2}\, \Big(W_{X,v,1}+\hat v_{\mu
_{q_i},j}^2\Big)\, =\, W_{X,v,1}$, $\ (J_{X,v})^{\mu _{q_i}-3}\,
\Big(W_{X,v,1}+\hat v_{\mu _{q_i},j}^2\Big) \, \neq  \, W_{X,v,1}$.
Moreover for $i=\tilde i +\tilde {\epsilon }+1, \ldots ,u$ and
$j\in \{ q_i-q_{i-1},\ldots ,1\}$ we have that $\hat v_{\mu
_{q_i},j}^{\mu _{q_i}}=v_{\mu _{q_i},j}^{\mu _{q_i}}\not \in
W_{X,v,1}$, since $v_{\mu _{q_i},j}^{\mu _{q_i}}$ is written in the
same row of the graph of $R_B$ as one of the elements of $\Delta
_B^{\circ }$; this means that $(J_{X,v})^{\mu _{q_i}-1}\, \Big(
W_{X,v,1}+\hat v_{\mu _{q_i},j}^1\Big) \neq W_{X,v,1}$, while
 $(J_{X,v})^{\mu _{q_i}}\, \Big( W_{X,v,1}+\hat v_{\mu
_{q_i},j}^1\Big) = W_{X,v,1}$ since $J^{\mu _{q_i}}\, v_{\mu
_{q_i},j}^1=0$.\hspace{4mm} $\square $ \vspace{2mm}\newline For
any $(X,v)\in \cSE^{\star } _B$ we consider a basis $\widehat
{\Delta }_{B,X,v}$ with the property expressed in lemma \ref{next'}
and we define a map from $\widehat{\Delta } _{B,X,v}$ to $\{
e_1,\ldots ,e_{n-\omega_1} \} $ in the following way: if
$W_{X,v,1}+\hat{v}_{\mu_{q_i},j}^l$ is the $h-$th element of
$\widehat{\Delta }
 _{B,X,v}$ the image of $W_{X,v,1}+\hat{v}_{\mu _{q_i},j}^l$ is $e_h$. This map induces an isomorphism
 $p_{X,v}$ from ${\displaystyle  K^n/
 W_{X,v,1}}\, $ to $\, K^{n-\omega_1}$; we denote by $\pi _B$ the morphism from $\cSE^{\star } _B$ to $N(n-\omega_1,K)$
 defined by $$(X,v)\longmapsto p _{X,v} \,
\circ \, X_{v,1}\, \circ \, \big( p _{X,v}\big) ^{-1}\ .$$ Let $\cQ _B$ be
the open subset of $\cSE^{\star } _B$ of all the pairs $(X,v)$
such that the partition of $X$ is the maximum partition
$Q(B)=(\omega_1,\ldots ,\omega_z)$ which is associated to the
elements of $\cSE _B\, $ ($\cSN _B$).  Let $\widetilde {\cSE }
_{\widehat
 B}$ be the set of all the nilpotent endomorphisms of $K^{n-\omega _1}$ which are conjugated to elements of $\cSE
 _{\widehat B}$.
\begin{proposition} \label{end} The morphism $\pi _B$ has the following properties: \begin{itemize} \item[a)] $\cSE _{\widehat B}\, \subseteq \, \pi _B
\, \big( \cSE^{\star } _B\big)  \ $; \item[b)] $ \pi _B\, \big(
\cQ _B\big)\, \subseteq \, \widetilde {\cSE } _{\widehat
 B} \ $. \end{itemize}
\end{proposition}
\pf \begin{itemize} \item[a)]  We consider the pairs $(X,v)\in
\cSE^{\star }_B$ with the following properties: \begin{itemize} \item[$P_1)$] there are no arrows from elements of
$\Delta _B^{\circ }$ to elements of $\Delta _B-\Delta _B^{\circ }$; \item[$P_2)$]
$v=v_{\mu _{q_1},q_{1}}^1$.\end{itemize} Since $v_{\mu _{q_1},q_{1}}^1\in \Delta _B ^{\circ }$,
for these pairs $W_{X,v,1}= \langle \, \Delta _B ^{\circ }\, \rangle $;
hence, by lemma \ref{nextt}, $a(i,j,i',j',l')=0$ for all $i=1,\ldots ,\tilde i-1$,
$j=q_i-q_{i-1},\ldots ,1$ and $(i',j',l')\in H(i,j)$, that is
$\widehat {\Delta }_{B,X,v}= \{ W_{X,v,1}+v_{\mu _{q_i},j}^{l}\ |\ v_{\mu
_{q_i},j}^{l}\in \Delta _B-\Delta _B^{\circ }\} $ by Corollary
\ref{next'}. If $v_{\mu _{q_i},j}^l$ and $v_{\mu
_{q_{i'}},j'}^{l'}$ are respectively the $h-$th and the $h'-$th
vector of $\Delta _B-\Delta _B^{\circ }$ then the entry of
$p _{X,v} \, \circ \, X_{v,1}\, \circ \, \big( p _{X,v}\big) ^{-1}\, (e_h)\, $
with respect to $e_{h'}$ is the entry of $X\, \big( v_{\mu
_{q_i},j}^l\big) $ with respect to $v_{\mu _{q_{i'}},j'}^{l'}$. By the definition of
$R_B$, Proposition \ref{order} and Corollary \ref{graph} we get that
for all $\widehat X\in \cSE _{\widehat B}$ there
exists $(X,v)\in \cSE^{\star  } _B$ with the properties $P_1)$ and $P_2)$ such that $\pi _B
((X,v))=\widehat{X}$. \item[b)] Let $(X,v)\in \cQ _B$ and
let $(X',v')\in \cQ _B\, \cap \, \big( \cN _B \times K^n)$. By
Theorem \ref{indToep} and Proposition \ref{j2} the endomorphisms
$X_{v,1}$ and $(X')_{v',1}$ are conjugated. Since $(X')_{v',1}$ commutes with $
J_{X',v'}$ we have that $\pi _B ((X',v'))\in \cN _{\widehat B}$,
hence $\pi _B ((X,v))\in \widetilde {\cSE } _{\widehat
 B}$.\end{itemize} \hspace{4mm} $\square $\vspace{2mm}
\newline  Now we can prove Theorem \ref{eend} of section \ref{introduction}.\vspace{2mm} \newline {\bf Proof of Theorem
\ref{eend}}\vspace{2mm}\newline {\em $1^{st}$ step}:\
by Theorem \ref{indToep} the maximum
partition which is associated to the elements of $\widetilde {\cSE
}_{\widehat B}\, $ ($\cSE _{\widehat B}$) is $(\widehat{\omega}_1,
\widehat {\omega}_2, \ldots ,\widehat{\omega}_{\widehat z})$.\vspace{2mm}\newline {\em $2^{nd}$ step}:\
 by Propositions \ref{f1} and by a) of Proposition \ref{end} there exists a not empty
open subset of $\cSE _B$ such that if $X$ belongs to it the
partition of $X_{v,1}$ is greater than or equal to $(\widehat{\omega}_1, \widehat {\omega}_2,
\ldots ,\widehat{\omega}_{\widehat z})$; hence by Proposition
\ref{j2} there exists a not empty
open subset of $\cSE _B$ such that if $X$ belongs to it the partition of $X$ is greater than or equal to
$(\omega _1,\widehat{\omega}_1, \widehat {\omega}_2,
\ldots ,\widehat{\omega}_{\widehat z})$.\vspace{2mm}\newline {\em $3^{rd}$ step}:\  by b) of Proposition \ref{end} the maximum partition
which is associated to the elements
of $\cSE _B$ is less than or equal to $(\omega_1,\widehat{\omega}_1, \widehat{\omega}_2,
\ldots ,\widehat{\omega}_{\widehat z})$. \vspace{2mm}\newline Hence the maximum partition which is associated to the elements
of $\cSE _B$ is $(\omega_1,\widehat{\omega}_1, \widehat{\omega}_2,
\ldots ,\widehat{\omega}_{\widehat z})$, which by Theorem
\ref{indToep} is equivalent to the claim. \hspace{4mm} $\square $\vspace{2mm} \newline
{\bf Another proof of Theorem \ref{eend}} We observe that all the proves of this Section could be repeated by substituting the variety $\cSE _B$ with the variety $\cSN _B$, getting a proof of Theorem \ref{eend} which is independent from Theorem \ref{indToep}. \hspace{4mm} $\square $.
\section{Other results on the orbits
intersecting ${\cN }_B$}\label{fifth}  In \cite{BasI} one proved a relation
between $Q(B)$ and the Hilbert functions of the algebras $K[A,J]$
for $A\in \cN _B$. After this result Toma\v{z} Ko\v{s}ir and
Polona Oblak in \cite{Kos} (2008) proved that a generic pair of
commuting nilpotent matrices generates a Gorenstein algebra and
then obtained as a consequence the following result.
\begin{theorem}\label{KO} For any partition $B$ the partition
$Q(B)$ has decreasing parts differing by at least 2, hence the map
$Q$ is idempotent. \end{theorem} Theorem \ref{eend} leads to the
following proof of Theorem \ref{KO}. \vspace{2mm} \newline {\bf
Another proof of Theorem \ref{KO} } We can prove the claim by
induction on $n$, hence we can assume that $\widehat{\omega
}_{i}-\widehat{\omega }_{i+1}\geq 2$ for $i=1,\ldots ,\widehat
z-1$. By the definition of $\omega _1$ and $\widehat{\omega }_1$
we get that $\omega _1-\widehat {\omega }_1\geq 2$, hence by
Theorem \ref{eend} we get that $\omega _i-\omega _{i+1}\geq 2$ for
$i=1,\ldots ,z-1$. Then by Corollary \ref{corR} we get that
$Q(Q(B))=Q(B)$. \hspace{4mm} $\square $ \vspace{2mm}
\newline Polona Oblak recently published a paper on properties of
the nilpotent orbits which intersect $\cN _B$ for some special
types of $B$ (see \cite{Obll}), among which we cite the following
result.
\begin{theorem} \label{Obll} If $n>3$ then $\cN _B$ intersects all the nilpotent
orbits if and only if $J^2=0$.
\end{theorem}
The map $Q$ was investigated by D.I. Panyushev in the more general
context of Lie algebras. In fact if $\mathfrak{g}$ is a semisimple
Lie algebra over an algebraically closed field $K$ such that char
$K=0$ and $\cN (\mathfrak g)$ is the nilpotent cone of
$\mathfrak{g}$ then $\cN (\mathfrak g)$ is irreducible (see
\cite{Ko}).  Let $e\in \cN (\mathfrak{ g})$ and let
$\mathfrak{z}_{\mathfrak g}(e)$ be the centralizer of $e$. Let $\{
e,h,f\} $ be an $\mathfrak{sl}_2$-triple and let
$\mathfrak{g}={\displaystyle \bigoplus _{i\in \Z }\;
\mathfrak{g}(i)}$ be the corresponding $\Z-$ grading of
$\mathfrak{g}$.  Let $G$ be the adjoint group of $\mathfrak{g}$.
In \cite{Pan} D.I. Panyushev defined an element $e$ to be
self-large if $G\cdot e\ \cap \ (\mathfrak{z}_{\mathfrak g}(e)\,
\cap \, \cN (\mathfrak g))$ is open (dense) in
$\mathfrak{z}_{\mathfrak g}(e)\, \cap \, \cN (\mathfrak g)$ and he
 proved the following result.
\begin{theorem}\label{PP} The element $e\in \cN $ is "self-large"
iff $\mathfrak{z}_{\mathfrak g}(e)\ \cap \ \mathfrak g (0)$ is
toral and $\mathfrak{z}_{\mathfrak g}(e)\ \cap \ \mathfrak g
(1)=\{ 0\} $. \vspace{4mm}\end{theorem}

{\small{\em Authors'
e-mail address:
robasili@alice.it}
\end{document}